\documentclass [12pt] {article}
\title { Quantum equations for knots}
\author{Thomas Fiedler  \& Appendix by Butian Zhang}
\usepackage{amsmath, amsfonts, amssymb, latexsym, graphicx}
\usepackage{caption}
\usepackage{subcaption}
\usepackage{tikz}
\usepackage{tikz-cd} 
\usepackage{breqn}
\usepackage{url}
\usepackage{hyperref}

\begin{document}
\newtheorem{proposition}{Proposition}
\newtheorem{theorem}{Theorem}
\newtheorem{lemma}{Lemma}
\newtheorem{corollary}{Corollary}
\newtheorem{example}{Example}
\newtheorem{remark}{Remark}
\newtheorem{definition}{Definition}
\newtheorem{question}{Question}
\newtheorem{conjecture}{Conjecture}
\newtheorem{observation}{Observation}
\maketitle

\begin{abstract}

This paper contains linear systems of equations which can distinguish knots without knot invariants.

Let $M_n$ be the topological moduli space of all n-component string links and such that a fixed projection into the plane is an immersion. If a string link is the product of some string link diagram $T$ and the parallel n-cable of a framed long knot diagram $D$, then there is a canonical arc $push$ in $M_n$, defined by pushing $T$ through the n-cable of $D$. In this paper we apply the combinatorial 1-cocycles from the HOMFLYPT and Kauffman polynomials in $M_n$ with values in the corresponding skein modules to this canonical arc in $M_n$. Some of the 1-cocycles  lead to linear systems of equations in the skein modules, for each couple of diagrams $D$ and $D'$. If the system has no solution in the Laurent polynomials then $D$ and $D'$ represent different knots. 

{\em We give first examples where we distinguish knots without any knot invariants. 
In particular, we distinguish the knot $9_{42}$ from its mirror image with equations coming from the HOMFLYPT polynomial. 
Notice that the knot $9_{42}$ and its mirror image share the same HOMFLYPT polynomial.
}

On the other hand, each solution of the system gives rather fine information about any regular isotopy which connects $D$ with $D'$.

\end{abstract}

\tableofcontents

\section{Introduction}

In the last 40 years many new knot invariants were discovered, in particular the well known quantum knot invariants, see e.g. \cite{K}. All these invariants are defined by using a single  diagram of a knot. {\em We can go one step further and look for invariants of 1-parameter families of knot diagrams and which connect two given diagrams.}

Moduli spaces of knots (and not only with a fixed projection into the plan which is an immersion) were extensively studied in  \cite {Bu},  \cite {Bu2}, \cite {Gr}, \cite {H}, \cite {HK}. Vassiliev has introduced the finite type cohomology classes of the moduli space of long knots \cite {V1}, \cite {V2}. These cohomology classes can be represented by combinatorial formulas. In particular, in dimension 0 the finite type invariants can be represented by the Gauss diagram formulas which were introduced in  \cite {PV}. Combinatorial formulas for some classes in dimension 1 (which are conjectured to be of finite type)  were obtained in  \cite {T}, \cite {S}, \cite {KK}, \cite {Mor1}, \cite {Mor2}, \cite {Mor3}, \cite {F2}, \cite {F3}.
\vspace{0,4 cm}

Our tool in this paper are also  combinatorial 1-cocycles in $M_n$, but which are based on quantum polynomials. In Chapter 6 of the monograph  \cite{F3} we have constructed a series of 1-cocycles in $M_n$ by using the HOMFLYPT and the Kauffman polynomials.

The present paper is a continuation of the Chapter 6 and we apply some of these 1-cocycles not to loops in $M_n$ but to the arc $push$ and we consider their values in the corresponding skein modules. 
\vspace{0,4 cm}

It is well known that knot types of classical knots and of long knots are in 1-1 correspondence. Moreover, two long knots which are isotopic, which have the same writhe and the same Whitney index, are regularly isotopic, i.e. there exist a combinatorial isotopy without Reidemeister I moves. We replace the diagram $D$ of the long knot now by its parallel n-cable, which we denote by $nD$. Consequently, two diagrams $D$ and $D'$ with the same writhe and the same Whitney index represent the same knot in $\mathbb{R}^3$ if and only if $nD$ and $nD'$ are regularly isotopic string links. Here we consider the projection $pr: \mathbb{R}^3  \to \mathbb{R}^2$ given by $(x,y,z) \to (x,y)$, where the long knots coincides with the $x$-axes outside a compact set. For practical reasons of calculation we will consider mainly the cases $n=2$ and $n=3$.

The machinery for the construction of 1-cocycles in $M_n$ was laid in two monographs  \cite{F2} and  \cite{F3}, see also \cite{FK}. Roughly speaking, one has to solve the positive global tetrahedron equation, the cube equations and the commutativity relations. The interested reader should consult \cite{F2}. In  Sections 2.2 and 2.3 we will just remind the definitions of the 1-cocycles. The proofs that they are indeed 1-cocycles are contained in Chapter 6 of \cite{F3}.
\vspace{0,4 cm}

The points at $x=-\infty$ are numbered by the increasing $y$-coordinate and we declare the point $1$ as the distinguished point $\infty$. Abstractly, we close the oriented string link $nD$ now to an oriented circle by the cyclic permutation $\sigma_1 \sigma_2...\sigma_{n-1}$ of the points at $x=\infty$. Each crossing c of $nD$ has of course a local type, positive or negative, but it has now also a global type, denoted by $\partial c = A$, where $A$ is the subset of the points $\{1,2,..,n\}$ which are contained in the arc of the circle $nD$ from the over-cross of c to the under-cross of c. 

All our combinatorial 1-cocycles are constructed by a discrete integration in the following way: in a generic regular isotopy a finite number of Reidemeister III and Reidemeister II moves occur.  We associate an element in the
 
corresponding skein module to each Reidemeister III move and to each  Reidemeister II move, and we sum up the contributions in the regular isotopy.

\begin{figure}
\centering
\includegraphics{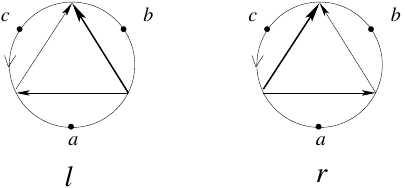}
\caption{\label{globtricross}  The six global types of triple crossings, e.g. $r_a$ means the type on the right with $a=1=\infty$.}  
\end{figure}

To each Reidemeister move of type III corresponds a diagram with a {\em triple 
crossing} $p$: three branches of the knot (the highest, middle and lowest with respect to the projection $pr$. A small perturbation of the triple crossing leads to an ordinary diagram with three crossings near $pr(p)$.
 We name the corresponding branches of the knot the lowest, middle and highest branch and the corresponding crossings $ml$ (lowest with middle), $hm$ (middle with highest) and $d$ (for distinguished). For better visualization we draw the crossing $d$ always with a thicker arrow. There are exactly six  global types of triple crossings with respect to the point $1=\infty$, shown in Fig.~\ref{globtricross} . The two new crossings from a Reidemeister II move are also called distinguished and denoted by $d$ too.

Let $\gamma(t), t \in [0,1],$ be a regular isotopy of $D$ to $D'$. It induces a regular isotopy from $nD$ to $nD'$, which we denote also by $\gamma(t)$. 

We had introduced the {\em trace graph} $TG(\gamma)$ as an immersed graph in the plan $(x,t)$ in \cite {FK}, see also  \cite {F2}. We follow a crossing $c_t$ in the isotopy and we orient the corresponding arc in $(x,t)$ in such a way that the local mapping degree for the projection into the $t$-axes is +1 if the crossing $c_t$ is positive and is -1 if it is negative. We glue the arcs together in Reidemeister II moves but we do not care about the Reidemeister III moves, i.e. either we consider the three crossings in the move as virtual or we consider a resolution by separating the branches of the triple point. The connected components of the resulting arcs are called {\em cobordism components of crossings} and are denoted by $c_t$.

The meridian of the stratum in the moduli space of a diagram with a quadruple crossing  consists of a sequence of eight R III moves. They come in couples corresponding to the remaining branch is on the left or the right of the triple crossing. The tetrahedron equations requires that the contributions of these eight moves cancel out together, see e.g. \cite{KaKo} for a thorough study of the local tetrahedron equation.

We name the four local branches in the quadruple crossing from the lowest to the highest (with respect to the $z$-coordinate for the projection into the plane) by 1 up to 4, and we name the unique crossing between two branches by the corresponding couple of numbers as well as the triple crossings by the corresponding triple of numbers.

{\em The heart of the matter is the fact that we had found  two very special kinds of solutions of the global tetrahedron equation, i.e. the contributions of the triple crossings depend on their global type too.}
\vspace{0,4 cm}

In the solutions of the first kind, the two moves $234$ cancel out together as well as the two moves $123$. The remaining moves $124$ and $134$ cancel only out all four together. The crucial point is that they share the distinguished crossing $d=14$ on the same cobordism component of crossings. Therefore we can restrict the 1-cocycle on each component of the cobordism of crossings for all those homotopies of the isotopy were surgeries (Morse modifications of index 1) of the cobordism do not occur. We achieve this by considering only very special loops and arcs in the moduli space $M_n$ together with very special homotopies, see Section 2.4. But this is enough in order to get invariants of cobordisms of crossings, see Section 4. 
\vspace{0,4 cm}

In the solutions of the second kind we use in addition some weight functions of finite type. This time only the two moves $234$ cancel out together. We can no longer restrict the 1-cocycles to cobordism components of crossings, because the crossing $d=14$ for the moves $124$ and $134$ is not the crossing $d=13$ of the moves $123$. However, there is still a splitting of our 1-cocycles by using the global types $A$ of crossings. The values of the 1-cocycles on the arc $push$ depend now very strongly on the string link $T$ and the type $A$. 

{\em For some choices of $T$ and $A$ this value is a knot invariant of $D$.} This is explained in Section 4 too.

{\em For other choices of $T$ and $A$ its value leads to a system of linear equations in the skein module. The variables of this system depend only on the couple of diagrams $D$ and $D'$ and are invariants of any regular isotopy (if there is any) which connects the diagrams $D$ and $D'$. For each basis element of the skein module the solution of the system (if there is any) is a rational function of two variables with rational coefficients. If this solution is not an integer Laurent polynomial of two variables then the diagrams $D$ and $D'$ represent different knots. If there is such a solution then it is an invariant of any regular isotopy (if there is any) which connects $D$ to $D'$, independent of the relative homotopy class of this regular  isotopy in $M_n$. }This is explained in Section 5.
\vspace{0,4 cm}

In Section 3 we give a first application of a version of our simplest 1-cocycle. There are two basic local types of triple crossings of oriented links: we call them {\em braid-like} (the only sort which occurs for braids) and {\em star-like} (which never occure for braids),  compare  Fig.~\ref{coorient}. It is well known, that star-like triple crossings can be avoided in each  1-parameter family of diagrams (i.e. arcs and loops in $M_n$). This follows essentially from the cube equations in \cite {F2}, see also \cite {P}.

Using our 1-cocycle $\tilde Y_c$ we show by an example that this is no longer the case for 2-parameter families: {\em there is a contractible loop in $M_1$ which does not contain star-like triple crossings, but each contracting 2-disc in $M_1$ contains necessarily star-like triple crossings.}

\begin{remark}
In fact, we had found our first 1-cocycle $Y_c$ already in 2008 and Dror Bar-Natan had written a computer program in order to calculate its value in the skein module $H_2(z,v)$  on an arc, called $scan$, see \cite{BN}. Unfortunately, this arc corresponds just to pushing $2D$ half through a positive curl on one branch. It is therefore not surprising, that the result was just the element in the skein module $H_2(z,v)$ represented by $2D$ multiplied with the HOMFLYPT polynomial of a positive curl, i.e. by $v$. (Because we are only interested in diagrams up to regular isotopy, we do not normalize the quantum invariants as usual by multiplication with $v^{-w(D)}$.) In this paper we push $nD$ through another fixed n-cable $nK$ and we observe, that this gives already an invariant of $nD$ for each fixed crossing of $nK$. We change then our point of view and push now a fixed string link $T$ through $nD$ and we refine  considerably the 1-cocycle in order to get invariants of $nD$ and to get our system of quantum equations!

\end{remark}

\section{The 1-cocycles}

\subsection{Generalities}

Let $\Sigma$ denote the discriminant of non-generic diagrams in $M_n$ together with its natural stratification, compare \cite{F2}.

\begin{figure}
\centering
\includegraphics{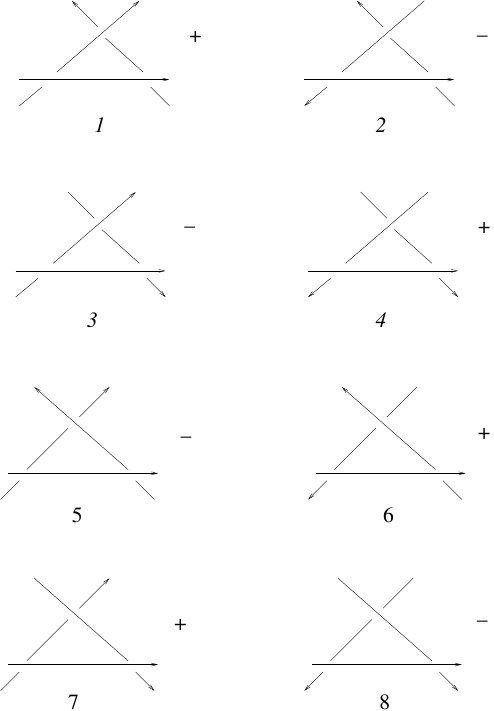}
\caption{\label{loctricross} Local types of a triple crossing}  
\end{figure}

For oriented knots there are eight local types of R III moves, compare Fig.~\ref{loctricross} and four local types of R II moves. The sign in the figure indicates the side of the discriminant $\Sigma_{tri}^{(1)}$ for the global type $r$, for the global type $l$ the sign is the opposite. Different local types of Reidemeister moves come together in strata of $\Sigma^{(2)}$ of codimension 2 in the moduli space $M_n$, compare Chapters 1 and 2 in  \cite {F2}. 

\begin{definition}

 Smoothing an ordinary crossing $q$ of a diagram $nD$ with respect to the orientation splits the (abstract) closure of $nD$ into two oriented and ordered circles. We call $D^+_q$ the component which goes from the over-cross to the under-cross at $q$ and by $D^-_q$ the remaining component, compare Fig.~\ref{splitq}.

\begin{figure}
\centering
\includegraphics{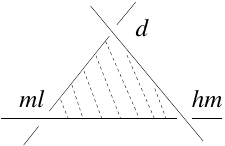}
\caption{\label{names} The names of the crossings in a R III-move}  
\end{figure}

\begin{figure}
\centering
\includegraphics{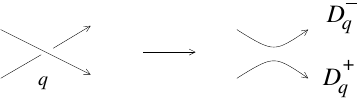}
\caption{\label{splitq} Two ordered knot diagrams associated to a crossing $q$}  
\end{figure}

\end{definition}

A {\em Gauss diagram of\/} a knot $nD$ is an oriented circle with oriented chords and with the points at infinity $\{1=\infty, 2, ...\}$. The chords correspond to the crossings of the knot diagram and are always oriented from the under-cross to the over-cross (here we use the orientation of the $z$-axes). Moreover, each chord (or arrow) has a sign, denoted by $w(c)$, which corresponds to the usual sign of the crossing $c$. The circle of a Gauss diagram in the plane is always equipped with the counter-clockwise orientation.

A {\em Gauss diagram formula\/} of degree $k$ is an expression assigned to the
diagram of a knot $nD$, which is of the following form:
\begin{displaymath}
\sum_{configurations}{\textrm{function( writhes of the crossings )}}
\end{displaymath}
where the sum is taken over all possible choices of $k$ (unordered)
different crossings in the diagram such that the chords arising
from these crossings in the diagram of $nD$ build a given sub-diagram
with respect to the points at infinity, but without fixing the signs on the arrows. The marked sub-diagrams are called
{\em configurations\/}, compare e.g. \cite{PV}, \cite{F2}.
The function, which is often the product of the writhes of the crossings in the configuration, is called the {\em weight}. In this paper only Gauss diagram formulas of degree $k=1$ will occur.

A Reidemeister III move corresponds to a triangle in the Gauss diagram. The {\em global type of a Reidemeister III move} was shown in Fig.~\ref{globtricross}.   We have to indicate whether the arrow $ml$ in the triangle goes to the left, denoted by $l$, or it goes to the right, denoted by $r$, and in which arc is the point $1=\infty$.

\begin{figure}
\centering
\includegraphics{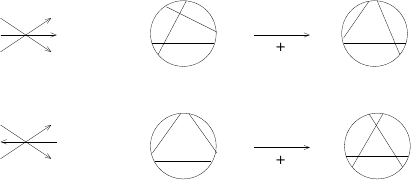}
\caption{\label{coorient} The coorientation for Reidemeister III-moves, braid-like and star-like}  
\end{figure}

\begin{definition}
The {\em coorientation} for a Reidemeister III move  is the direction from two intersection points of the corresponding three arrows to one intersection point and of no intersection point of the three arrows to three intersection points, compare 
Fig.~\ref{coorient}. (One can see in the cube equations for $\Sigma^{(2)}_{trans-self}$ that the two coorientations for triple crossings fit together for the strata of $\Sigma^{(1)}_{tri}$ which come together in $\Sigma^{(2)}_{trans-self}$, compare \cite{F2}.) Evidently, our coorientation is completely determined by the corresponding planar curves and therefore we can draw just chords instead of arrows in Fig.~\ref{coorient}.  We call the side of the complement of $\Sigma^{(1)}$ in $M_n$ into which points the coorientation, the {\em positive side} of  $\Sigma^{(1)}$.

The coorientation for Reidemeister II moves is the direction from no crossings to the diagram with two new crossing. 

\end{definition}

Each transverse intersection point $p$ of an oriented generic arc in $M$ with $\Sigma^{(1)}$ has now an intersection index $+1$ or $-1$, called $sign(p)$, by comparing the orientation of the arc with the coorientation of $\Sigma^{(1)}$.

\vspace{0,4 cm}

For the convenience of the reader we remind now our strategy developed in \cite{F2}.

We study the relations of the local types of triple crossings in what we call the {\em cube equations}.

\begin{figure}
\centering
\includegraphics{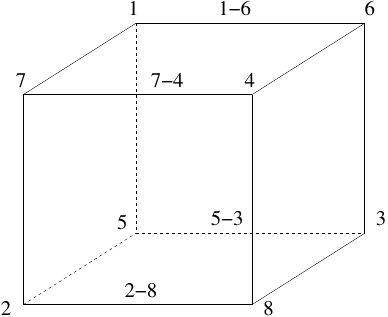}
\caption{\label{gamma} The graph $\Gamma$}  
\end{figure}

Triple crossings come together in points of $\Sigma^{(2)}_{trans-self}$, i.e. an auto-tangency with in addition a transverse branch \cite{F2}. But one easily sees that the global type of the triple crossings  is always preserved. We make now a graph $\Gamma$ for each global type of a triple crossing in the following way: the vertices correspond to the different local types of triple crossings. We connect two vertices by an edge if and only if the corresponding strata of triple crossings are adjacent to a stratum of $\Sigma^{(2)}_{trans-self}$. We have shown that the resulting graph is the 1-skeleton of the 3-dimensional cube $I^3$, see Fig.~\ref{gamma}.
In particular, it is connected. The edges of the graph $\Gamma = skl_1(I^3)$ correspond to the types of strata in $\Sigma^{(2)}_{trans-self}$. The solution of the positive tetrahedron equation tells us what is the contribution to the 1-cocycle of a positive triple crossing (i.e. all three involved crossings are positive) . The meridians of the strata from $\Sigma^{(2)}_{trans-self}$ give equations which allow us to determine the contributions of all other types of triple crossings. However, a global phenomenon occurs: each loop in $\Gamma$ could give an additional equation. Evidently, it suffices to consider the loops which are the boundaries of the 2-faces from $skl_2(I^3)$. We call all the equations which come from the meridians of $\Sigma^{(2)}_{trans-self}$ and from the loops in $\Gamma = skl_1(I^3)$ the {\em cube equations}. (Notice that a loop in $\Gamma$ is more general than a loop in $M_n$. For a loop in $\Gamma$ we come back to the same local type of a triple crossing but not necessarily to the same whole diagram of the knot.)
\vspace{0,4 cm}

The 2-sphere, which is the boundary of the cube, has an additional structure. The six edges which connect the braid-like triple crossings, i.e. types 1, 3, 4, 5, 7, 8, form an {\em equator} of the 2-sphere. The star-like types 2 and 6 form the {\em poles} of the 2-sphere and there are six edges which connect the poles with the equator.
\vspace{0,4cm}

Our strategy was the following: first we find a solution of the positive global tetrahedron equation, which gives us the contribution of the local type 1. We solve then the cube equations for the different local types of triple crossings and auto-tangecies. 

The commutation relations, which correspond to a meridian of the 

transverse intersection of two strata of codimension one, are automatically satisfied. Indeed, our partial smoothings of triple crossings and auto-tangecies as well as the weights of degree 1 are invariant under a simultaneous R III or R II move in another place of the diagram. We have shown in \cite{F2} that the resulting 1-cochain is then a 1-cocycle in $M_n$, because it satisfies automatically all other tetrahedron equations, which was shown by using strata of codimension three in the discriminant $\Sigma \subset M_n$.
\vspace{0,4cm}

It is amazing, that we have found solutions of the cube equations with integer coefficients in the case of the HOMFLYPT polynomial. However, we could find such solutions in the case of the Kauffman polynomial only on the edges of the equator. The edges to the poles force us to reduce the coefficients mod 2. We have no explication for this phenomenon. A guess would be, that this is perhaps related to 1-parameter families in representation theory. Indeed, it is well known that the HOMFLYPT polynomial for knots (and their cables) is related to the fundamental representation of $SU(n)$ and  the Kauffman polynomial is related to the fundamental representation of $SO(n)$. But $SU(n)$ is simply connected and the (stable) fundamental group of $SO(n)$ is $\mathbb{Z}/2\mathbb{Z}$.

\subsection{The case of the HOMFLYPT polynomial: $Y_c$ and $R_{reg}^{(1)}$}

\begin{figure}
\centering
\includegraphics{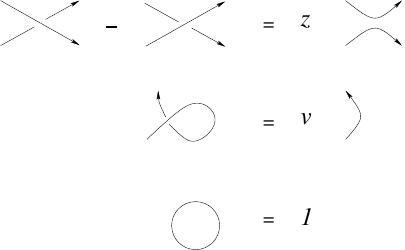}
\caption{\label{HOMFLYPT}  skein relations for the HOMFLYPT polynomial $P$}  
\end{figure}

The HOMFLYPT polynomial  polynomial is defined by the skein relations shown in Fig.~\ref{HOMFLYPT}. 

The skein module of n component string links can be identified with the Hecke algebra $H_n(z,v)$ associated to the braid group $B_n$, compare e.g. \cite{K}. We consider the usual basis ($id, \sigma_1$) for $n=2$ and  ($id, \sigma_1, \sigma_2, \sigma_1 \sigma_2, \sigma_2 \sigma_1, \sigma_1 \sigma_2 \sigma_1=\sigma_2 \sigma_1 \sigma_2$) for $n=3$.

\begin{definition}
An ordinary crossing $q$ in a diagram $nD$ is of {\em type 1} if $\infty \in D^+_q$ and is of {\em type 0} otherwise and we write $[q]=1$, respectively $[q]=0$.
\end{definition}

We consider a braid-like triple crossing $p$ locally like an element in $H_3(z,v)$ and we replace it by the {\em partial smoothing}  

$III(p)=(\sigma_2(p) - \sigma_1(p))w(d)w(hm)w(ml)$ if and only if  $[d]=0$.

If $p$ is a star-like triple crossing then its partial smoothing are the same as  shown in Fig.~\ref{smooth}.

If $p$ is an auto-tangency with opposite direction of the branches and with $[d]=0$ then its partial smoothing $II_+(p)$ is defined in Fig.~\ref{sts}.

Let $nK$ be the parallel n-cable of a fixed long knot diagram $K$. We want to define invariants of the parallel n-cable $nD$ of a long knot diagram $D$.

\begin{definition}
Let $c$ be a fixed crossing of $nK$ with $[c]=0$. We push $nD$ under the over-branch of $c$ and we call this arc $push(D,c)$ in $M_n$. Then the 1-cochain $Y_c$ with values in $H_n(z,v)$ is defined by 
\vspace{0,4cm}

$Y_c=\sum_{p,r_a,r_b,d=c} sign(p) III(p)-\sum_{p,l_b,d=c} sign(p) III(p)-\sum_{p,d=c} sign(p)II_+(p)$.
\vspace{0,4cm}

We define another 1-cochain without the contributions of auto-tangencies and only with the contributions of braid-like triple crossings by 
\vspace{0,4cm}

$\tilde Y_c=\sum_{p=braid-like,r_a,r_b,d=c} sign(p) III(p)-\sum_{p=braid-like,l_b,d=c} sign(p) III(p)$

\end{definition}

\begin{theorem}
$Y_c(push(D,c))$ in $H_n(z,v)$ is an invariant of $nD$ for the cobordism component $c_t$ of $c$ in $push(D,c)$.

Moreover, if $push(D,c)$ does not contain star-like triple crossings, then $\tilde Y_c(push(D,c))$ is already an invariant of the cobordism component $c_t$ in $push(D,c) \subset M_n \setminus \Sigma^{(1)}_{star-like}$.

\end{theorem}

{\em Proof.}

In \cite{F2}, Proposition 6.1, we haven proven that $Y_c$ and hence also $\tilde Y_c$ satisfy the positive global tetrahedron equations. As already explained in the Introduction, we can restrict the solution to the cobordism component of $c_t= d=14$.

It follows from an easy examination of the figures in \cite{F3}, that $\tilde Y_c$ and hence $Y_c$ satisfies the cube equations on the braid-like equator:

For the global type $r$: Fig. 6.52, Fig. 6.53, Fig. 6.54, Fig. 6.55, Fig. 6.56, Fig. 6.64.

For the global type $l$: Fig. 6.28, Fig. 6.34, Fig. 6.36, Fig. 6.38, Fig. 6.42, Fig. 6.45. 

Simultaneous Reidemeister moves at other places in the diagram do not change the contribution of $p$ at all.

Consequently, $\tilde Y_c$ is a 1-cocycle for each cobordism component $c_t$. Moreover, one easily sees that if we push $nD$ under the over-cross of a crossing $c$ in $nK$ then in all Reidemeister moves in $push(D,c)$ the crossing $c_t$ is the distinguished crossing $d$ of the move. 
\vspace{0,4cm}

The edges to the poles are exactly the edges for which the branches in the auto-tangencies have opposite tangent direction. They have to contribute now too, in order to give a solution of the cube equations.

Inspecting the figures Fig. 6.57, Fig. 6.58, Fig. 6.59, Fig. 6.60, Fig. 6.62, Fig. 6.63 for the global type $r$ and Fig. 6.30, Fig. 6.32, Fig. 6.40, Fig. 6.46, Fig. 6.48, Fig. 6.51 for the global type $l$, one easily sees that  $Y_c$ satisfies all cube equations now. 

Again, they are well defined on the cobordism components $c_t=d$. Indeed, the contributions of the two auto-tangecies in an edge do not cancel out together if and only if the third branch passes between the two branches of the auto-tangency. But exactly in this case the crossing $c_t=d$ from the auto-tangencies is also the crossing $c_t=d$ in the two triple crossings in the edge.
\vspace{0,4cm}

The contributions from the auto-tangencies are only needed on the edges to the star-like triple crossings. Consequently, $\tilde Y_c$ on $push(D,c)$ is an invariant in 
$push(D,c) \subset M_n \setminus \Sigma^{(1)}_{star-like}$.

$\Box$

\vspace{0,4cm}

We remind now the definition of the more sophisticated 1-cocycle $R_{reg}^{(1)}$, which uses six from the eight strata in the tetrahedron equation.

Let in all generality $T(p)$ be a diagram of a string link (with the abstract closure to an oriented circle  by a given identification of the end points) and which has  a triple crossing or a self-tangency $p$  and let $d$ be the distinguished crossing for $p$. In the case of a self-tangency we identify the two distinguished crossings. {\em We assume that the distinguished crossing $d$ is of type 0.}

\begin{definition}
A crossing $q$ of $T(p)$ is called a {\em f-crossing} if $q$ is of type 1 and  the under cross of $q$ is in the oriented arc from $\infty$ to the over cross of $d$ in $T(p)$.
\end{definition}

\begin{figure}
\centering
\includegraphics{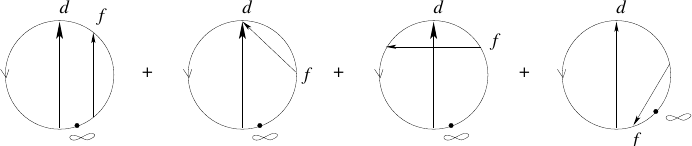}
\caption{\label{foot}  the f-crossings}  
\end{figure}

\begin{figure}
\centering
\includegraphics{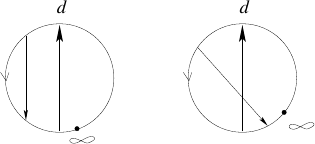}
\caption{\label{notfoot} crossings of type 1 which are not f-crossings}  
\end{figure}

\begin{figure}
\centering
\includegraphics{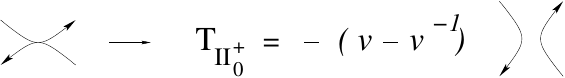}
\caption{\label{sts}  partial smoothing of a self-tangency $II_0^+$}  
\end{figure}

In other words, the "foot" of $q$ in the Gauss diagram is in the sub arc of the circle which goes from $\infty$ to the head of $d$ and the head of $q$ is in the sub arc which goes from the foot of $q$ to $\infty$.   We illustrate this in Fig.~\ref{foot} and Fig.~\ref{notfoot}.  Notice that we make essential use of the fact that no crossing can ever move over the points at infinity.

\begin{definition}
Let $q$ be a f-crossing. The {\em global type} $A=\partial q$ of $q$ is defined by $\partial q = D^+_q \cap \partial T$. 
\end{definition}

The global type $\partial q$ is a cyclically ordered subset of $\partial T$  which contains $\infty$ (compare the Introduction).

\begin{definition}
Let a type $A$ be fixed. The {\em linear weight $W_A(p)$} is  the Gauss diagram invariant shown in  Fig.~\ref{foot} too. In other words, it is the sum of the signs of all f-crossings of global type $A$ with respect to $p$ in a given diagram $T(p)$. Notice that if $p$ is a self-tangency then the degenerate second configuration in  Fig.~\ref{foot} can not appear.
\end{definition}

\begin{lemma}
The weight $W_A(p)$ is an  isotopy invariant for each Reidemeister move  of type II or III, i.e. $W_A(p)$ is invariant under any isotopy of the rest of the string link outside a small neighborhood of $p$.
\end{lemma}

{\em Proof.} This is obvious. One has only to observe that the two new crossings from a Reidemeister move of type II are either both not f-crossings or are both f-crossings. In the latter case they have the same type $A$ but they have different signs. 

$\Box$

 It is clear that $W_A(p)$ is an invariant of order 1 for the couple $(T,p)$ for each fixed type $A$.
\vspace{0,4cm}

Let $\gamma$ be an oriented generic arc in $M_n$ which intersects $\Sigma^{(1)}$ only in positive triple crossings and let $A$ be a fixed global type.

\begin{definition}
The evaluation of the 1-cochain $R^{(1)}_{reg}$ with global type $A$ on $\gamma$  is defined by \vspace{0,2 cm}

 $R^{(1)}_{reg}(A)(\gamma)=\sum_{p,r_a} sign(p) \sigma_2\sigma_1(p) -\sum_{p,l_c} sign(p) \sigma_2\sigma_1(p)$

$+ \sum_{p,r_a,r_b}sing(p)zW_A(p)(\sigma_2(p)-\sigma_1(p))$

$-\sum_{p,l_b}sing(p)zW_A(p)(\sigma_2(p)-\sigma_1(p))$     \vspace{0,2 cm}

where the first sum is over all triple crossings of the global types shown in Fig.~\ref{globtricross} (i.e. the types $l_c$ and $r_a$) and such that $\partial  (hm)=A$.  The second sum is over all triple crossings $p$ which have a distinguished crossing $d$ of type 0 (i.e. the types $r_a$, $r_b$ and $l_b$).

\end{definition}

(Remember that e.g. $\sigma_2\sigma_1(p)$ denotes the element in the HOMFLYPT skein module $H_n(z,v)$ which is obtained by replacing the triple crossing $p$ in the string link  $T$ by the partial smoothing $\sigma_2\sigma_1$.)

Notice that a triple crossing of the first type $r_a$ in Fig.~\ref{globtricross} could also contribute to the second sum in $R^{(1)}_{reg}(A)(\gamma)$ but a triple crossing of the second type $l_c$ in Fig.~\ref{globtricross} can't because $d$ is of type 1.

$R^{(1)}_{reg}(A)$ is a solution of the positive global tetrahedron equations. It has to be extended to a 1-cocycle by a solution of the cube equations.

For a triple crossing $p$ of a given local type from $\{1,...,8\}$ we denote by $T_{r_a}(type)(p)=T_{l_c}(type)(p)$ and $T(type)(p)$ the corresponding partial smoothing of $p$ in the skein module $H_n(z,v)$. We will determine these partial smoothing from the cube equations.

\begin{definition}
We distinguish four types of self-tangencies, denoted by 

$II_0^+, II_0^-, II_1^+, II_1^-$. Here "0" or "1" is the type of the distinguished crossing $d$ and "+" stands for opposite tangent directions of the two branches and "-" stands for the same tangent direction.
\end{definition}

\begin{definition}
The partial smoothing $T_{II^+_0}(p)$ of a self-tangency with opposite tangent direction is defined in Fig.~\ref{sts}.
The partial smoothing $T_{II^-_0}(p)$ of a self-tangency with equal tangent direction is defined as $-zP_T$ (where $P_T$ is the element represented by the original string link  $T$ in the HOMFLYPT skein module $H_n$).
\end{definition}

Notice that the distinguished crossing $d$ has to be of type 0 in both cases.

\begin{definition}
Let $p$ be a self-tangency. Its contribution to $R^{(1)}_{reg}(A)$ is defined by

 $R^{(1)}_{reg}(A)=sign(p) W_A(p)T_{II^+_0}(p)$ respectively $sign(p) W_A(p)T_{II^-_0}(p)$.\vspace{0,2 cm}

\end{definition}

\begin{figure}
\centering
\includegraphics{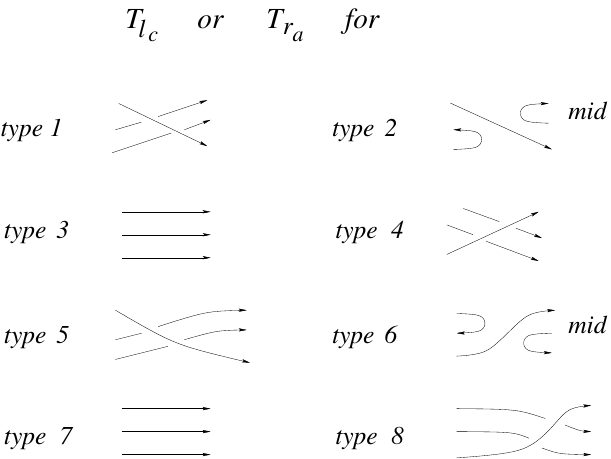}
\caption{\label{typesmooth}  the partial smoothings $T_{l_c}$ and $T_{r_a}$}  
\end{figure}

\begin{figure}
\centering
\includegraphics{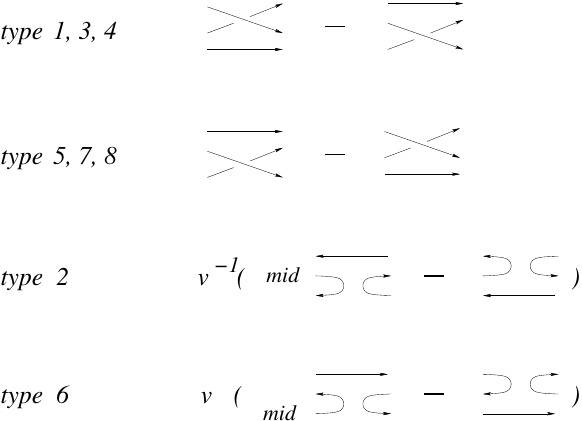}
\caption{\label{smooth} the partial smoothings $T$ (type)}  
\end{figure}

\begin{definition}
The partial smoothing for the local and global types of triple crossings are given in Fig.~\ref{typesmooth} and  Fig.~\ref{smooth} (remember that {\em mid} denotes the ingoing middle branch for a star-like triple crossing). 
\end{definition}

\begin{definition}
Let $s$ be a generic oriented arc in $M_n$ (with the given  abstract closure of the string link to an oriented circle). Let $A \subset \partial T$ be a given global type.  

The 1-cochain $R^{(1)}_{reg}(A)$ is defined by \vspace{0,2 cm}

$R^{(1)}_{reg}(A)(s ) = +\sum_{p, r_a }sign(p)T_{r_a}(type)(p) -\sum_{p, l_c} sign(p)T_{l_c} (type)(p) $\vspace{0,2 cm}

 $+\sum_{p, r_a,r_b }sign(p)zW_A(p)T(type)(p) $\vspace{0,2 cm}

 $-\sum_{p, l_b }sign(p)zW_A(p)T(type)(p) $\vspace{0,2 cm}

 $+\sum_{p \in II^+_0 } sign(p)W_A(p)T_{II^+_0}(p)$\vspace{0,2 cm} 

 $+\sum_{p \in  II^-_0 } sign(p)W_A(p)T_{II^-_0}(p)$\vspace{0,2 cm} 

Here all weights $W_A(p)$ are defined only over the f-crossings $f$ with $\partial f= A$ and in the first two sums (i.e. for $T_{l_c}$ and $T_{r_a}$) we require that $\partial (hm)=A$ for the triple crossings.

Notice that the stratum $r_a$ contributes in two different ways.
\end{definition}

\begin{theorem}

The 1-cochain $R^{(1)}_{reg}(A)$ with values in $H_n(z,v)$ is a 1-cocycle in $M_n$ for each global type $A$.

\end{theorem}

This is a part of Theorem 6.2 in \cite{F3}.

\subsection{The case of the 2-variable Kauffman polynomial: $Z_c$ and $R_{F,reg}^{(1)}$}

\begin{figure}
\centering
\includegraphics{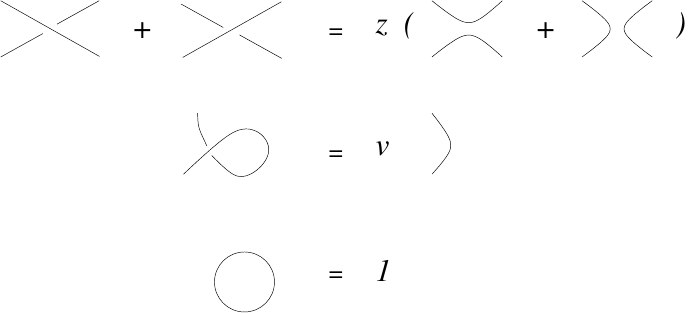}
\caption{\label{basisKauf} skein relations for the Kauffman polynomial $F$}  
\end{figure}

The skein relations for the 2-variable Kauffman invariant are shown in Fig.~\ref{basisKauf}, see e.g. \cite{K}.

\begin{figure}
\centering
\includegraphics{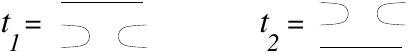}
\caption{\label{3kauf} new generators of the Kauffman skein module of 3-component string links}  
\end{figure}

Let $t_1$ and $t_2$ be the new generators of the Kauffman skein module of 3-component string links as shown in Fig.~\ref{3kauf}. We call the corresponding new generator of the Kauffman  skein module of 2-component string links also by $t_1$, i.e. without the third branch.

The Kauffman skein module of n-component string links can be identified with the Birman-Murakami-Wenzl algebra $C_n(z,v)$ associated to the braid group $B_n$, compare e.g. \cite{K}. We consider the usual basis ($id, \sigma_1, t_1$) for $n=2$ and  ($id, \sigma_1, \sigma_2, \sigma_1 \sigma_2, \sigma_2 \sigma_1, \sigma_1 \sigma_2 \sigma_1=\sigma_2 \sigma_1 \sigma_2,t_1,t_2,t_1t_2,t_2t_1$) for $n=3$.

The types, the signs, the weights, the global type $A$ in the Kauffman case are exactly the same as in the HOMFLYPT case. Only the partial smoothing are different.

We consider a triple crossing $p$ locally like an element in $C_3(z,v)$ and an auto-tangency locally like an element in $C_2(z,v)$ and we replace them by the {\em partial smoothings}. 

It turns out that this time $(\sigma^2_2(p)-\sigma^2_1(p))$ is a solution of the positive global tetrahedron equation, see Lemma 6.16 in \cite{F3}.
\vspace{0,4cm}

But surprisingly there is second solution of the positive global tetrahedron equation, namely $(\sigma_1^{-1}t_2+\sigma_2t_1-t_1\sigma_2^{-1}-t_2\sigma_1)$, see Lemma 6.18 in \cite{F3}. We have solved the cube equations for this solution. Here is the result. We leave the verification to the reader, by using the same figures as in the HOMFLYPT case.
\vspace{0,4cm}

\begin{definition}

type 1 = $\sigma_1\sigma_2\sigma_1$: $III(p)=\sigma_1^{-1}t_2+\sigma_2t_1-t_1\sigma_2^{-1}-t_2\sigma_1$
\vspace{0,4cm}

type 3 = $\sigma_1\sigma_2^{-1}\sigma_1^{-1}$: $III(p)=(z+v^{-1})t_1-(z+v^{-1})t_2 +t_2\sigma_1-\sigma_2t_1$
\vspace{0,4cm}

type 4 = $\sigma_1^{-1}\sigma_2^{-1}\sigma_1$: $III(p)=(z+v^{-1})t_1-(z+v^{-1})t_2 +\sigma_1t_2-t_1\sigma_2$
\vspace{0,4cm}

type 5 = $\sigma_1^{-1}\sigma_2\sigma_1$: $III(p)=vt_1-vt_2+\sigma_2t_1-t_2\sigma_1$
\vspace{0,4cm}

type 7 = $\sigma_1\sigma_2\sigma_1^{-1}$: $III(p)=vt_1-vt_2+
t_1\sigma_2-\sigma_1t_2$
\vspace{0,4cm}

type 8 = $\sigma_1^{-1}\sigma_2^{-1}\sigma_1^{-1}$: $III(p)=-\sigma_1^{-1}t_2-\sigma_2t_1+t_1\sigma_2^{-1}+t_2\sigma_1$
\vspace{0,4cm}

type 2 = $(w(d)=+1,w(ml)=w(hm)=-1)$: $\sigma_1-\sigma_2+t_2\sigma_1-\sigma_2t_1$
\vspace{0,4cm}

type 6 = $(w(d)=-1,w(ml)=w(hm)=+1)$: $\sigma_1^{-1}-\sigma_2^{-1}+\sigma_1^{-1}t_2-t_1\sigma_2^{-1}$
\vspace{0,4cm}

$II_0^+$: $II_0^+(p)=t_1$
\vspace{0,4cm}

 $II_0^-$: $II_0^-(p)=t_1(v-v^{-1})/z$
\end{definition}

\begin{definition}
Let $c$ be a fixed crossing of a fixed diagram $nK$ with $[c]=0$. As for $Y_c$  in the HOMFLYPT case we push $nD$ under the over-branch of $c$ in $nK$, and we call this arc as usual $push(D,c)$ in $M_n$. Then the 1-cochain $Z_c$ with values in $C_n(z,v)$ is defined by 
\vspace{0,4cm}

$Z_c=\sum_{p,r_a,r_b,d=c} sign(p) III(p)-\sum_{p,l_b,d=c} sign(p) III(p)$
\vspace{0,2cm}

$+\sum_{p,d=c} sign(p)II_+(p)+\sum_{p,d=c} sign(p)II_-(p)$.
\vspace{0,4cm}

\end{definition}

We have the analogue of Theorem 1 for the Kauffman polynomial too, but only with coefficients in $\mathbb{Z}/2\mathbb{Z}$.

\begin{theorem}
$Z_c(push(Dc))$ in $C_n(z,v)$ with coefficients  in  $\mathbb{Z}/2\mathbb{Z}$ is an invariant of $nD$ for the cobordism component $c_t$ of $c$  in $push(D,c) \subset M_n$.

\end{theorem}
\vspace{0,4cm}

Let us consider now  solutions  of the global positive tetrahedron equation in Kauffman's case and which use six from the eight strata in the tetrahedron equation.

Notice that we could rewrite equivalently the  solution of the global positive tetrahedron equation in the case of the HOMFLYPT polynomial  in the following way : \vspace{0,2 cm}

$R^{(1)}_{reg}(A)(\gamma)=\sum_{p,r_a} sign(p) \sigma_2\sigma_1(p) -\sum_{p,l_c} sign(p) \sigma_2\sigma_1(p)$
\vspace{0,2cm}

$+ \sum_{p,r_a,r_b}sing(p)W_A(p)(\sigma^2_2(p)-\sigma^2_1(p))-\sum_{p,l_b}sing(p)W_A(p)(\sigma^2_2(p)-\sigma^2_1(p))$  \vspace{0,2 cm}

where the first two  sums are only over the global types $r_a$ and $l_c$ of triple crossings and with $\partial (hm)=A$.

{\em Amazingly, exactly the same formula is also a solution in the case of the Kauffman polynomial!}

However, we haven't solved the cube equations in this case. They are very complicated. Instead we will consider another 1-cochain, which is related to the unexpected solution $(\sigma_1^{-1}t_2+\sigma_2t_1-t_1\sigma_2^{-1}-t_2\sigma_1)$.
We have the following proposition, compare Proposition 6.10 in \cite{F3}.

\begin{proposition} The 1-cochain \vspace{0,2 cm}

$R^{(1)}_{F,reg}(\gamma)(A)=\sum_{p,r_a} sign(p) t_1t_2(p)-\sum_{p,l_c} sign(p) t_1t_2(p)$
\vspace{0,2cm}

$ + \sum_{p,r_a,r_b}sing(p)W_A(p)z(-\sigma_1^{-1}t_2-\sigma_2t_1+t_1\sigma_2^{-1}+t_2\sigma_1)$
\vspace{0,2cm}

$-\sum_{p,l_b}sing(p)W_A(p)z(-\sigma_1^{-1}t_2-\sigma_2t_1+t_1\sigma_2^{-1}+t_2\sigma_1)$ 
\vspace{0,2 cm}

is also a solution of the positive global tetrahedron equation. Here the first two sums are over all triple crossings of the global types $l_c$ and $r_a$ and such that $\partial  (hm)=A$.  The remaining sums are over all triple crossings $p$ which have a distinguished crossing $d$ of type 0 (i.e. the types $r_a$, $r_b$ and $l_b$).

\end{proposition}

\begin{figure}
\centering
\includegraphics{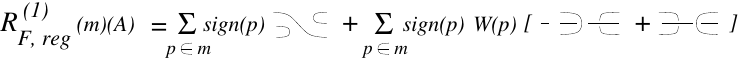}
\caption{\label{solkauf}  surprising solution of the positive tetrahedron equation in Kauffman's case}  
\end{figure}

(compare also Fig.~\ref{solkauf}, where we give it in a shorter form for locally defined signs)

\begin{figure}
\centering
\includegraphics{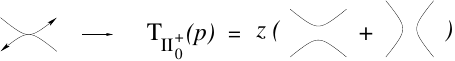}
\caption{\label{Ksts} partial smoothing of a $RII$ move with opposite tangent directions}  
\end{figure}

\begin{figure}
\centering
\includegraphics{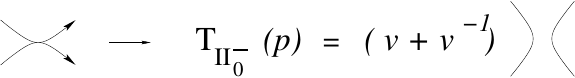}
\caption{\label{Keq} partial smoothing of a $RII$ move with the same tangent direction}  
\end{figure}

\begin{definition}
For $R^{(1)}_{F, reg}$:

The partial smoothing $T_{II^+_0}(p)$ of a self-tangency with opposite tangent direction and $d$ of type 0 is defined in Fig.~\ref{Ksts}.

The partial smoothing $T_{II^-_0}(p)$ of a self-tangency with equal tangent direction and $d$ of type 0 is defined in Fig.~\ref{Keq}.

\end{definition}

\begin{figure}
\centering
\includegraphics{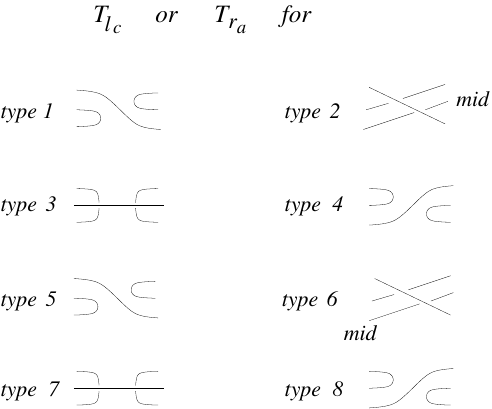}
\caption{\label{Ktypesmooth} the partial smoothings $T_{l_c}$ and $T_{r_a}$ in the
Kauffman case}  
\end{figure}

\begin{figure}
\centering
\includegraphics{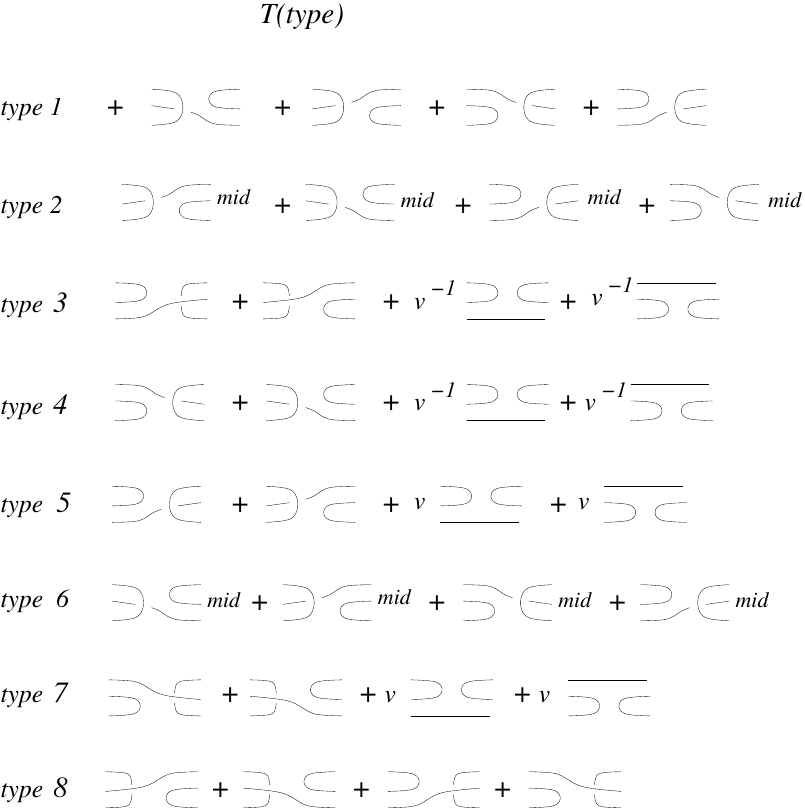}
\caption{\label{Ksmooth} the partial smoothings $T(type)$ in the
Kauffman case}  
\end{figure}

We will solve now the cube equations.

\begin{definition}
The partial smoothings for the local and global types of triple crossings in Kauffman's case are given in Fig.~\ref{Ktypesmooth} and  Fig.~\ref{Ksmooth} (remember that {\em mid} denotes the ingoing middle branch for a star-like triple crossing). 

\end{definition}

\begin{definition}
Let $s$ be a generic oriented arc in $M_n$ (with a fixed abstract closure of the string link to an oriented circle). Let $A$ be a given global type.  The 
1-cochain $R^{(1)}_{F, reg}(A)$ is defined by \vspace{0,2 cm}

$R^{(1)}_{F,reg}(A)(s ) = \sum_{p,r_a} sign(p)T_{r_a} (type)(p) -\sum_{p,l_c}sign(p)T_{l_c}(type)(p) $\vspace{0,2 cm}

 $+\sum_{p,r_a,r_b }sign(p)zW_A(p)T(type)(p) 
 -\sum_{p, l_b }sign(p)zW_A(p)T(type)(p) $\vspace{0,2 cm}

 $+\sum_{p \in  II^+_0 } sign(p)W_A(p)T_{II^+_0}(p)
+\sum_{p \in  II^-_0 } sign(p)W_A(p)T_{II^-_0}(p).$\vspace{0,2 cm}

Here all weights $W_A(p)$ are defined only over the f-crossings $f$ with $\partial f= A$ and in the first two sums (i.e. for $T_{l_c}$ and $T_{r_a}$) we require that $\partial (hm)=A$ for the triple crossings.

\end{definition}

\begin{theorem}
The 1-cochain $R^{(1)}_{F, reg}(A)$ with coefficients in $\mathbb{Z}/2\mathbb{Z}$ and with values in $C_n(z,v)$ is a 1-cocycle in $M_n$ for each global type $A$.

\end{theorem}

This is a part of Proposition 6.16 in \cite{F3}.

\subsection{Nice arcs, loops and homotopies of loops in the moduli space $M_n$}

Let $T$ be an (oriented) n-component string link diagram and $nD$ be the parallel n-cable of a  long knot diagram. We denote by $TnD$ their usual product, i.e. we glue $nD$ on the right to $T$. Let $\gamma$ be a regular isotopy of $D$ to another long knot diagram $D'$. We denote the induced regular isotopy of $TnD$ to $TnD'$ by $T\gamma$.

The arc $push(T,nD)$ starts in $TnD$ in $M_n$. We push  $T$ monotonously through $nD$ up to $nDT$.
\vspace{0,4cm}

We glue the following arcs to a loop in $M_n$, called $push(T,D,D')$, which depends only on $T$, $D$ and $D'$:
\vspace{0,4cm}

$push(T,D,D')=push(T,nD)$  $(\gamma T)$  $(-push(T,nD'))$  $(-T \gamma )$.
\vspace{0,4cm}

The important point is, that this loop is {\em nicely} contractible in $M_n$. Indeed, $push(T,nD)$ and $(T\gamma)$ commute and hence we have a 1-parameter family of loops, called  $push(T,D,D')_t$, by shortening $(T\gamma)$ with respect to the parameter $t$ of the arc $\gamma$. We end up with a loop which is just $push(T,nD)$  followed by $-push(T,nD)$ which we shorten  in itself to the point $TnD$.   

It follows that all our 1-cocycles  $R^{(1)}_{reg}(A)$ and $R^{(1)}_{F,reg}(A)$ vanish on the loop $push(T,D,D')$. 

\begin{remark}
The arc $push(T,nD)$ can be completed to a non-trivial loop in $M_n$ by adding $push(nD,T)$ if and only if the string link $T$ is itself the parallel n-cable of a long knot diagram.

\end{remark}

\section{A first application of $\tilde Y_c$}

We consider the figure-eight knot as a closed 3-braid $\sigma_2 \sigma_1^{-1}\sigma_2 \sigma_1^{-1}$ and we open up the third strand to make it a long knot. This is our diagram $D$. As $D'$ we take the mirror image of $D$. The figure-eight knot is amphichiral and $D$ and $D'$ share the same writhe and Whitney index. Consequently, there is a regular isotopy $\gamma$ in $M_1$ which connects $D$ with $D'$ and which does not contain star-like triple crossings, compare the Introduction.

\begin{figure}
\centering
\includegraphics{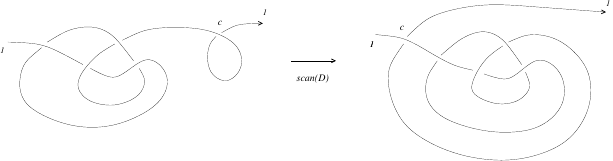}
\caption{\label{figureequations}  the scan-arc for $D$}  
\end{figure}

We add a small positive curl on the right of $D$ and of $D'$ on the branch which goes to infinity and we push $D$, respectively $D'$, just along the over-cross of the crossing in the curl. This arc in $M_1$ is called $scan(D)$, respectively $scan(D')$, compare Fig.~\ref{figureequations}. 

We consider the analogue for $scan$ of our contractible loop $push(T,D,D')$. Evidently, the contributions to $\tilde Y_c$ in $\gamma$ and $-\gamma$ cancel out together, because in both cases the curl adds just a factor $v$ for each contribution. It follows that if in the contraction of the loop there do not appear any star-like triple crossings, then $\tilde Y_c(scan(D))=\tilde Y_c(scan(D'))$.
\vspace{0,4cm}

In $scan(D)$ there are exactly four triple crossings. We give them with sign, global type, locale type:

$-, l_b,type 7$

$+,l_b,type 1$

$+,r_a, type 7$

$+,l_b,type 1$

An easy calculation gives now the following contributions to $\tilde Y_c$:

$-(v^{-1}\delta P(H^+) -v^{-2}P(3_1^+))$

$-(v^{-1}P(4_1) - v \delta P(H^-))$

$-(v - v^{-1} \delta P(H^+))$

$-(v^{-1} - v \delta P(H^-))$

Here as usual, $\delta = (v-v^{-1})/z$, $3_1^+$ the positive trefoil, $4_1$ the figure-eight knot, $H^+$ and $H^-$ the positive and the negative Hopf link.

Consequently, $\tilde Y_c(scan(D))= 2v(v-v^{-1})^2z^{-2}$ plus terms of higher degrees in $z$.
\vspace{0,4cm}

In $scan(D')$ there are also exactly four triple crossings:

$-,r_a, type 1$

$+,r_a, type 7$

$+,l_b, type 1$

$+,r_a, type 7$

An easy calculation gives now the following contributions:

$-(v^{-1}(\delta P(H^+) -zv^{-1}P(H^+))-v^{-1})$

$-(v - v^{-1} \delta P(H^+))$

$-(v^{-1}-v\delta P(H^-))$

$-(vP(4_1)-v^{-1} \delta P(H^+))$

Consequently, $\tilde Y_c(scan(D'))= (v+v^{-1})(v-v^{-1})^2z^{-2}$ plus terms of higher degrees in $z$.
\vspace{0,4cm}

It follows that star-like triple crossings have to appear in each contracting homotopy of the loop in $M_n$.

\section{The invariants from 1-cocycles applied to the arc $push$}

For our 1-cocycles $R^{(1)}_{ reg}(A)$ and $R^{(1)}_{F, reg}(A)$ we can no longer consider the cobordism components $c_t$ and on the other hand the global types of the triple crossings become essential. We know that the 1-cocycles vanish on each loop $push(T,D,D')$.
Consequently, the value of the 1-cocycle on $push(T,nD)$ (which we can calculate) is a knot invariant if and only if its values on $(T\gamma)$ and $(-\gamma T)$ cancel out together. We do not know anything about $\gamma$ besides the fact that $(T\gamma)$ and $(-\gamma T)$ contain exactly the same local Reidemeister moves, but with opposite signs. Hence, the pairs of moves will cancel out together if and only if they have also the same global type and the same weights $W_A(p)$. 

{\em One easily sees, that they share the same global type (including $\partial(hm)=A$) if and only if the string link $T$ induces the trivial permutation on its end points, and that they share the same weights $W_A(p)$ for a given type $A$ if no 1-crossing of $T$ is of the global type $A$.}

We call these two conditions on $T$ together the {\em $T$-neutral condition}.
\vspace{0,4cm}

For any string link $T$ we denote by $[T]$ the element represented by $T$ in the corresponding skein module, and we develop it in the standard basis (compare Sections 2.2 and 2.3): $[T]= T_{id} id+T_{\sigma_1} \sigma_1+...$, where all the coefficients $T_{id}$, $T_{\sigma_1}$,... are Laurent polynomials.

As well known, the multiplication in the skein modules $H_2(z,v)$ and $C_2(z,v)$ is commutative. For $n>2$ we will use only that each element commutes with $id$ and with itself and we consider only a very particular case, which generalizes the case for $n=2$.

The following theorem is an immediate consequence, because the contributions of $(T\gamma)$ and $(-\gamma T)$ cancel now out together.

\begin{theorem}
If a 2-component string link $T$ satisfies the $T$-neutral condition, then $R^{(1)}_{ reg}(A)(push(T,2D)) \in H_2(z,v)$ and $R^{(1)}_{F, reg}(A)(push(T,2D)) \in C_2(z,v)$ (with coefficients in $\mathbb{Z}/2\mathbb{Z}$ here) are invariants of $D$ up to regular isotopy of long knots.

For $n>2$, let $T$ be a n-component string link such that $[T]=T_{id}id+T_{\sigma_1}\sigma_1 \in H_n(z,v)$ (as always the case for $n=2$), respectively $[T]=T_{id}id+T_{\sigma_1}\sigma_1 +T_{t_1}t_1 \in C_n(z,v)$, and which satisfies the $T$-neutral condition (e.g. $T=\sigma_1^2$ and $A$ different from  $A=\{1=\infty,3,...,n\}$, which is the global type of the unique 1-crossing in $T=\sigma_1^2$).

Then $R^{(1)}_{ reg}(A)(push(T,nD))_{id}id+R^{(1)}_{ reg}(A)(push(T,nD))_{\sigma_1}\sigma_1$ 

and $R^{(1)}_{ F,reg}(A)(push(T,nD))_{id}id+R^{(1)}_{F, reg}(A)(push(T,nD))_{\sigma_1}\sigma_1$

$+R^{(1)}_{F, reg}(A)(push(T,nD))_{t_1}t_1$  (with coefficients in $\mathbb{Z}/2\mathbb{Z}$ here) are invariants of $D$ up to regular isotopy of long knots.

\end{theorem}

It is clear that the theorem could be generalized for more complicated string links $T$, but the formulation would also become more complicated.

\section{The quantum equations from 1-cocycles applied to the arc $push$}

We keep the condition that $T$ induces the trivial permutation on its end points, but we drop the condition on the weights $W_A(p)$. Consequently, the contributions of $(T\gamma)$ and $(-\gamma T)$ cancel no longer out together. However, we can express the sum of the values of our 1-cocycles on them by introducing some variables, which depend {\em only} on the existence of $\gamma$ and not on $T$. By developing the variables in the skein module and varying the string link $T$ we obtain linear systems of equations in the skein module. This will be our {\em quantum equations}.

The global type $A$ of the crossing $hm$ in a triple crossing of $\gamma$ does not depend on the position of $T$. Therefore we have only to control de changing of $W_A(p)$ with respect to the two positions of $T$. Remember that the crossing $d$ in the triple crossing or the auto-tangency has to be a 0-crossing.
\vspace{0,4cm}

{\em The HOMFLYPT case  for n=2}
\vspace{0,4cm}

\begin{figure}
\centering
\includegraphics{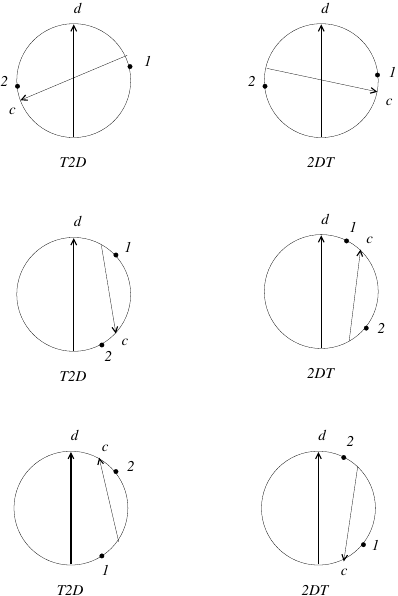}
\caption{\label{cd} the $d$ in the first two lines are of type $(1d)$}  
\end{figure}

For $T = T_1 =\sigma_1^2$ the unique 1-crossing $c$ in $T$ has the global type $A=\{1=\infty\}$ and hence we consider this $A$. We show the possible positions of $c$ with respect to $d$ in Fig.~\ref{cd}. In the first two cases $c$ is a f-crossing for $p$ in $T\gamma$ but not in $\gamma T$. In the third case $c$ is a f-crossing for both of them. (The arrow $c$ just moves on the circle to have its endpoints again near to the markings at infinity.) We say in the first two cases that $p$ is of {\em type} $(1d)$, i.e. following the orientation of the circle we come from $1=\infty$ first to the over-cross of $d$.

\begin{definition}

The variable $x \in H_2(z,v)$ for $\gamma$ is defined by 
\vspace{0,2 cm}

 $x=\sum_{p, r_a,r_b, (1d)}sign(p)zT(type)(p) $\vspace{0,2 cm}

 $-\sum_{p, l_b,(1d)}sign(p)zT(type)(p) $\vspace{0,2 cm}

 $+\sum_{p \in II^+_0,(1d) } sign(p)T_{II^+_0}(p)$\vspace{0,2 cm} 

 $+\sum_{p \in  II^-_0,(1d) } sign(p)T_{II^-_0}(p)$.\vspace{0,2 cm} 

Here we sum up over all Reidemeister moves $p$ of the type $(1d)$ in the regular isotopy of the parallel 2-cable induced by $\gamma$, compare Definitions 12 and 13.
\end{definition}

We decompose $x$ as usual in the skein module $x=x_{id}id+x_{\sigma_1} \sigma_1$, where $x_{id}$ and 
$x_{\sigma_1}$ are Laurent polynomials in $z$ and $v$.
\vspace{0,2 cm}

We denote in general $R^{(1)}_{ reg}(A)(push(T,nD))-R^{(1)}_{ reg}(A)(push(T,nD'))$ shortly by 
$\Delta(T)$, but which depends of course on $n$. This is the quantity which we can calculate.
\vspace{0,2 cm}

$R^{(1)}_{ reg}(A)(push(T,D,D'))=0$ and the changing of $W_A(p)$ for the moves $p$ of the type $(1d)$ implies now
\vspace{0,2 cm}

$\Delta(T)=Tx=xT \in H_2(z,v)$
\vspace{0,2 cm}

After developing in the skein module $H_2(z,v)$ we obtain our first system of quantum equations:
\vspace{0,4 cm}

$x_{id}+zx_{\sigma_1}=\Delta_{id}(T_1)$

$zx_{id}+(1+z^2)x_{\sigma_1}=\Delta_{\sigma_1}(T_1)$
\vspace{0,4 cm}

$det=1$ and consequently the system has a unique solution in the polynomials of the variables $z$ and $v$.

{\em The solution $x \in H_2(z,v)$ is an invariant of any regular isotopy which connects the diagrams $D$ and $D'$. }
\vspace{0,2 cm}

\begin{figure}
\centering
\includegraphics{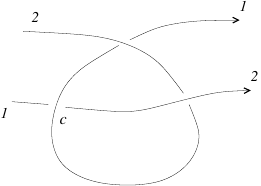}
\caption{\label{t2}  the auxiliary string link $T_2$}  
\end{figure}

Let us consider the string link $T_2$ shown in Fig.~\ref{t2}. The unique 1-crossing $c$ is negative and it is still of global type $A=\{1=\infty\}$. Consequently, this time
\vspace{0,2 cm}

$\Delta(T_2)=-T_2x=-xT_2 \in H_2(z,v)$.
\vspace{0,2 cm}

Developing in $H_2(z,v)$ we obtain another system for the same variables $x_{id}$ and $x_{\sigma_1}$
\vspace{0,4 cm}

$-vx_{id}+zv^{-1}x_{\sigma_1}=\Delta_{id}(T_2)$

$zv^{-1}x_{id}+(-v+z^2v^{-1})x_{\sigma_1}=\Delta_{\sigma_1}(T_2)$
\vspace{0,4 cm}

$det= (-1-v^{-2})z^2+v^2$ 
and hence the system has also a unique solution in the rational functions 
and which should be the {\em same} as the solution of the previous system, 
in particular, it should be in $\mathbb{Z}[v^{\pm 1}, z^{\pm 1}]$. 

{\em If this is not the case, then $D$ and $D'$ represent different knots!}
\vspace{0,2 cm}

The quantum equations can a priori distinguish knots, 
but in any case they give information about {\em any} regular isotopy connecting two given knot diagrams.
\vspace{0,4cm}

{\em The HOMFLYPT case for n=3}
\vspace{0,4cm}

The multiplication in the skein module $H_3(z,v)$ is no longer commutative.

We consider again $T=\sigma_1^2$. The unique 1-crossing $c$ in $T$ has the global type $A=\{1=\infty, 3\}$ now. Let us consider first the T-neutral case, i.e. $A=\{1=\infty\}$ or $A=\{1=\infty, 2\}$.

\begin{definition}

The new variable $y \in H_3(z,v)$ for $\gamma$ is defined by 
\vspace{0,2 cm}

$y=R^{(1)}_{ reg}(A)(\gamma)$, i.e. we evaluate the 1-cocycle $R^{(1)}_{ reg}(A)$ on the parallel 3-cable of $\gamma$, without $T$.

\end{definition}

$R^{(1)}_{ reg}(A)(push(T,D,D'))=0$ implies now 
\vspace{0,4 cm}

$\Delta(T)=Ty-yT \in H_3(z,v)$
\vspace{0,2 cm}

By developing in the skein module $H_3(z,v)$ we obtain the quantum equations
\vspace{0,4 cm}

$0=\Delta_{id}(T)$

$0=\Delta_{\sigma_1}(T)$

$y_{\sigma_1 \sigma_2}-y_{\sigma_2 \sigma_1}=z^{-1}\Delta_{\sigma_2}(T)$

$y_{ \sigma_2}+zy_{\sigma_1\sigma_2}-y_{\sigma_1\sigma_2 \sigma_1}=z^{-1}\Delta_{\sigma_1\sigma_2}(T)$

$-y_{ \sigma_2}-zy_{\sigma_2\sigma_1}+y_{\sigma_1\sigma_2 \sigma_1}=z^{-1}\Delta_{\sigma_2\sigma_1}(T)$

$-y_{\sigma_1 \sigma_2}+y_{\sigma_2 \sigma_1}=z^{-1}\Delta_{\sigma_1\sigma_2\sigma_1}(T)$
\vspace{0,4 cm}

This gives a system of four quantum equations with four variables. $det=0$ and if $\Delta_{\sigma_2}(T)+\Delta_{\sigma_1\sigma_2\sigma_1}(T) \not= 0$ then it has no solution, implying that $D$ and $D'$ represent different knots.
\vspace{0,4 cm}

Let us consider now the most interesting case $A=\{1=\infty, 3\}$.
One easily sees, that again $c$ is a f-crossing for $p$ in $T\gamma$ but not in $\gamma T$ if and only if $p$ is of type $(1d)$.

The variable $x \in H_3(z,v)$ is defined by exactly the same formula as in Definition 19, but of course now we sum up over all Reidemeister moves $p$ of the type $(1d)$ in the regular isotopy of the parallel 3-cable induced by $\gamma$.

$R^{(1)}_{ reg}(A)(push(T,D,D'))=0$ and the changing of $W_A(p)$ for the moves $p$ of the type $(1d)$ implies now
\vspace{0,2 cm}

$\Delta(T)=T(x+y)-yT \in H_3(z,v)$
\vspace{0,2 cm}

By developing in $H_3(z,v)$ we obtain a system of six equations with twelve variables. Hence we have to add more equations by varying the string link $T$, e.g. by taking $T_2$ for the first two strands as for $n=2$. The unique 1-crossing $c$ in $T_2$ is negative and has also the global type $A=\{1=\infty, 3\}$. But this should be best done with a computer program too.
\vspace{0,4cm}

{\em The Kauffman case for n=2}
\vspace{0,4cm}

For $T=\sigma_1^2$ the unique 1-crossing $c$ in $T$ has the global type $A=\{1=\infty\}$ and hence we consider this $A$.

The variable $x \in C_2(z,v)$ is defined by exactly the same formula as in Definition 19, but with the corresponding partial smoothing in Kauffman's case, compare Definitions 17 and 18.

$R^{(1)}_{F, reg}(A)(push(T,D,D'))=0$ and the changing of $W_A(p)$ for the moves $p$ of the type $(1d)$ implies now
\vspace{0,2 cm}

$\Delta(T)=Tx=xT \in C_2(z,v)$
\vspace{0,2 cm}

for $\Delta(T)$ defined by using of course $R^{(1)}_{F, reg}(A)$ instead of $R^{(1)}_{reg}(A)$.

After developing in the skein module $C_2(z,v)$ we obtain the following system of quantum equations with coefficients $mod$ 2:
\vspace{0,4 cm}

$(1+zv^{-1})x_{id}+zx_{\sigma_1}=\Delta_{id}(T)$

$zx_{id}+(1+z^2)x_{\sigma_1}=\Delta_{\sigma_1}(T)$

$(zv^{-2}+z^2v^{-1})x_{\sigma_1}+v^{-2}x_{t_1}=\Delta_{t_1}(T)$
\vspace{0,4 cm}

$det=v^{-2}+zv^{-3}+z^3v^{-3}$ and hence the system has a unique solution in the rational functions $\mathbb{Z}/2\mathbb{Z}(z,v)$.

Again, this solution is an invariant of any regular isotopy which connects the diagrams $D$ and $D'$. If they are not Laurent polynomials, then the diagrams represent different knots.

\section{Appendix by Butian Zhang}
In this appendix, we first give an example of the calculation of $R_{reg}^{1}$ of 2-cables
explaining how the quantum equations work to distinguish knots. 
Next, we give results on more examples. 

\subsection{How the quantum equations work ?}

Consider the following knot diagrams $D = 6_1$ and $D' = 6_2$. They have the same writhe $= 2$ and the same Whitney index $=0$. 

\begin{figure}[htbp]
  \centering
  \begin{subfigure}[t]{0.25\textwidth}
    \centering
    \includegraphics[height=1.1cm]{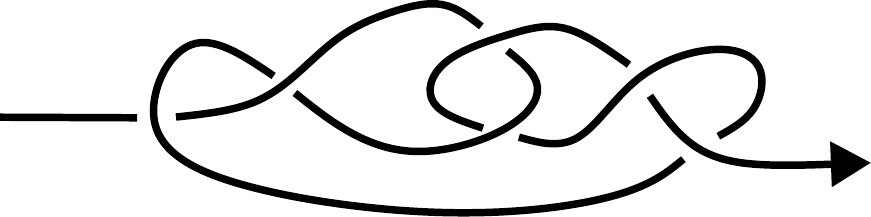}
    \caption{$D = 6_1$}
    \label{fig:6_1}
  \end{subfigure}
  \hspace{0.2\textwidth}
  \begin{subfigure}[t]{0.25\textwidth}
    \centering
    \includegraphics[height = 1.1cm]{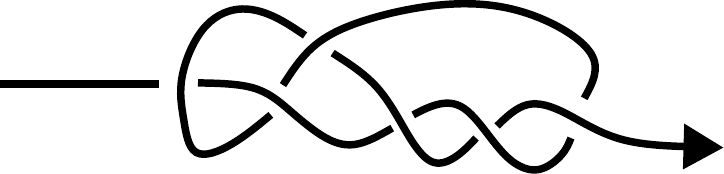}
    \caption{$D' = 6_2$}
    \label{fig:6_2}
  \end{subfigure}
  \caption{Two knot diagrams with the same writhe and Whitney index}
  \label{fig:knot_diagram_pair_6_1_6_2}
\end{figure}
We have the push arcs for the tangle $T_2$ on the parallel 2-cables of the two knot diagrams (blackboard framed).
\[
\begin{tikzcd}[column sep=3.5cm, row sep=1.5cm]
\includegraphics[width=5.5cm]{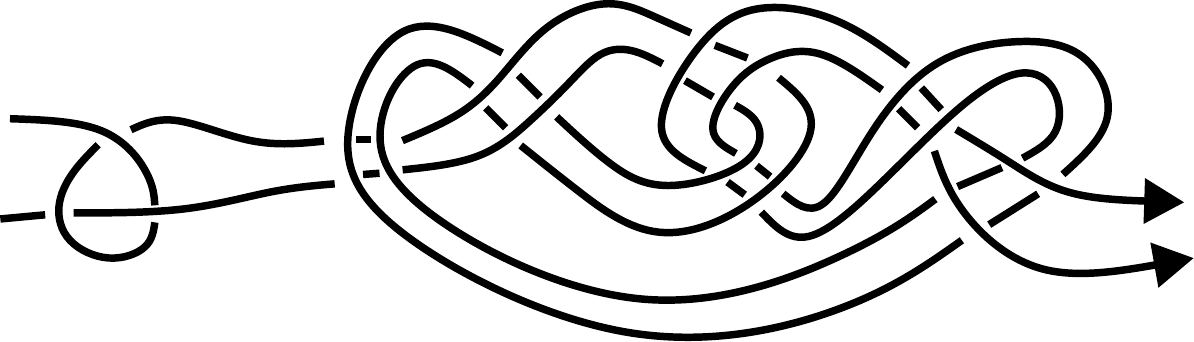}
  \arrow[r, "{\mathrm{push}(T_2, 2D)}"]
  \arrow[d, "\gamma"'] &
\includegraphics[width=5.5cm]{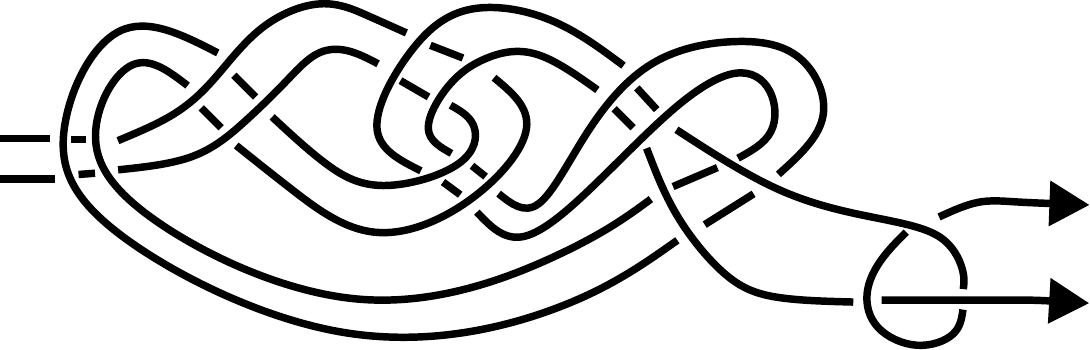}
  \arrow[d, "\gamma'"] \\
\includegraphics[width=5.5cm]{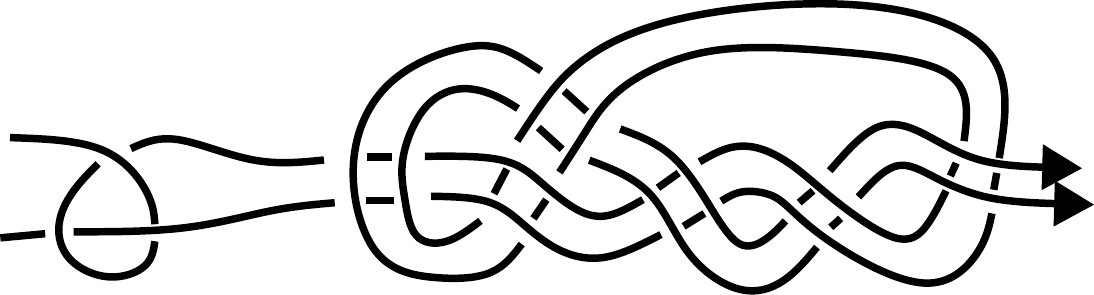}
  \arrow[r, "{\mathrm{push}(T_2, 2D')}"'] &
\includegraphics[width=5.5cm]{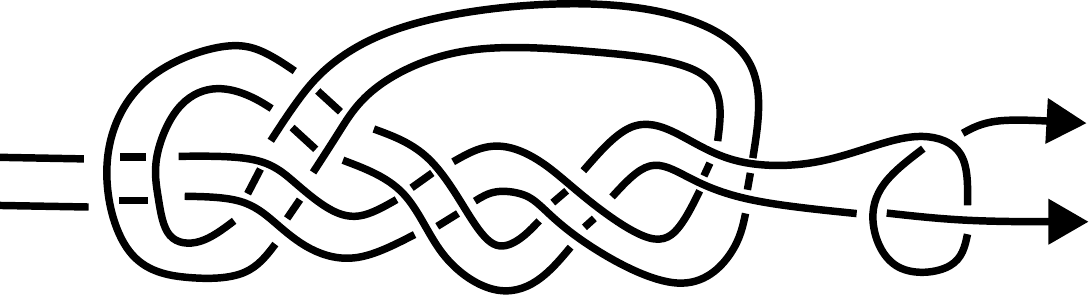}
\end{tikzcd}
\]

Using a computer program (see \cite{ZBT99}), we can get the following results.

\begin{dmath*}
R_{reg}^{1}(push(T_2, 2D)) = (3 v^{5} z^{2} + 4 v^{5} - 5 v^{3} z^{4} - 13 v^{3} z^{2} - 5 v^{3} + v z^{6} + 6 v z^{4} + 7 v z^{2} + 2 v - \frac{4 z^{8}}{v} - \frac{24 z^{6}}{v} - \frac{44 z^{4}}{v} - \frac{38 z^{2}}{v} - \frac{13}{v} - \frac{2 z^{10}}{v^{3}} - \frac{17 z^{8}}{v^{3}} - \frac{32 z^{6}}{v^{3}} + \frac{17 z^{4}}{v^{3}} + \frac{48 z^{2}}{v^{3}} + \frac{15}{v^{3}} + \frac{z^{12}}{v^{5}} + \frac{7 z^{10}}{v^{5}} + \frac{20 z^{8}}{v^{5}} + \frac{49 z^{6}}{v^{5}} + \frac{93 z^{4}}{v^{5}} + \frac{58 z^{2}}{v^{5}} + \frac{7}{v^{5}} + \frac{z^{12}}{v^{7}} + \frac{9 z^{10}}{v^{7}} + \frac{31 z^{8}}{v^{7}} + \frac{46 z^{6}}{v^{7}} + \frac{15 z^{4}}{v^{7}} - \frac{22 z^{2}}{v^{7}} - \frac{10}{v^{7}} - \frac{z^{8}}{v^{9}} - \frac{9 z^{6}}{v^{9}} - \frac{21 z^{4}}{v^{9}} - \frac{12 z^{2}}{v^{9}})
\includegraphics[height=1em]{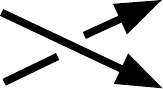} \linebreak
+ (2 v^{5} z + \frac{v^{5}}{z} + v^{3} z^{5} - v^{3} z^{3} - 6 v^{3} z - \frac{2 v^{3}}{z} - 2 v z^{5} - 8 v z^{3} - 8 v z - \frac{v}{z} + \frac{5 z^{9}}{v} + \frac{38 z^{7}}{v} + \frac{95 z^{5}}{v} + \frac{100 z^{3}}{v} + \frac{44 z}{v} + \frac{5}{v z} - \frac{13 z^{7}}{v^{3}} - \frac{68 z^{5}}{v^{3}} - \frac{96 z^{3}}{v^{3}} - \frac{39 z}{v^{3}} - \frac{3}{v^{3} z} + \frac{z^{11}}{v^{5}} - \frac{z^{9}}{v^{5}} - \frac{48 z^{7}}{v^{5}} - \frac{148 z^{5}}{v^{5}} - \frac{132 z^{3}}{v^{5}} - \frac{33 z}{v^{5}} - \frac{2}{v^{5} z} + \frac{z^{11}}{v^{7}} + \frac{4 z^{9}}{v^{7}} + \frac{z^{7}}{v^{7}} + \frac{4 z^{5}}{v^{7}} + \frac{37 z^{3}}{v^{7}} + \frac{28 z}{v^{7}} + \frac{2}{v^{7} z} + \frac{3 z^{7}}{v^{9}} + \frac{14 z^{5}}{v^{9}} + \frac{17 z^{3}}{v^{9}} + \frac{4 z}{v^{9}})
\includegraphics[height=1em]{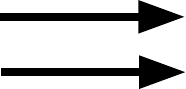}
\end{dmath*}

\begin{dmath*}
R_{reg}^{1}(push(T_2, 2D')) = (3 v^{5} z^{6} + 20 v^{5} z^{4} + 20 v^{5} z^{2} + 6 v^{5} - v^{3} z^{8} - 13 v^{3} z^{6} - 39 v^{3} z^{4} - 33 v^{3} z^{2} - 6 v^{3} - 5 v z^{12} - 50 v z^{10} - 172 v z^{8} - 251 v z^{6} - 213 v z^{4} - 119 v z^{2} - 24 v + \frac{2 z^{12}}{v} + \frac{31 z^{10}}{v} + \frac{171 z^{8}}{v} + \frac{424 z^{6}}{v} + \frac{489 z^{4}}{v} + \frac{251 z^{2}}{v} + \frac{49}{v} + \frac{5 z^{12}}{v^{3}} + \frac{42 z^{10}}{v^{3}} + \frac{109 z^{8}}{v^{3}} + \frac{69 z^{6}}{v^{3}} - \frac{64 z^{4}}{v^{3}} - \frac{87 z^{2}}{v^{3}} - \frac{27}{v^{3}} - \frac{4 z^{10}}{v^{5}} - \frac{31 z^{8}}{v^{5}} - \frac{74 z^{6}}{v^{5}} - \frac{59 z^{4}}{v^{5}} - \frac{11 z^{2}}{v^{5}} + \frac{2}{v^{5}})
\includegraphics[height=1em]{crossed_arrows.pdf} \linebreak
+ (2 v^{5} z^{5} + 17 v^{5} z^{3} + 13 v^{5} z + \frac{2 v^{5}}{z} - 4 v^{3} z^{9} - 40 v^{3} z^{7} - 132 v^{3} z^{5} - 151 v^{3} z^{3} - 65 v^{3} z - \frac{9 v^{3}}{z} + 5 v z^{13} + 58 v z^{11} + 255 v z^{9} + 555 v z^{7} + 687 v z^{5} + 474 v z^{3} + 153 v z + \frac{17 v}{z} - \frac{3 z^{13}}{v} - \frac{47 z^{11}}{v} - \frac{274 z^{9}}{v} - \frac{741 z^{7}}{v} - \frac{958 z^{5}}{v} - \frac{599 z^{3}}{v} - \frac{175 z}{v} - \frac{16}{v z} - \frac{4 z^{13}}{v^{3}} - \frac{38 z^{11}}{v^{3}} - \frac{120 z^{9}}{v^{3}} - \frac{122 z^{7}}{v^{3}} + \frac{47 z^{5}}{v^{3}} + \frac{141 z^{3}}{v^{3}} + \frac{62 z}{v^{3}} + \frac{6}{v^{3} z} + \frac{4 z^{11}}{v^{5}} + \frac{34 z^{9}}{v^{5}} + \frac{96 z^{7}}{v^{5}} + \frac{106 z^{5}}{v^{5}} + \frac{45 z^{3}}{v^{5}} + \frac{5 z}{v^{5}})
\includegraphics[height=1em]{parallel_arrows.pdf}
\end{dmath*}

so we have 

\begin{align*}
\Delta(T_2) &= R_{reg}^{1}(push(T_2, 2D)) - R_{reg}^{1}(push(T_2, 2D')) \\
&=\Delta_{id}(T_2) \includegraphics[height=1em]{parallel_arrows.pdf}
+ \Delta_{\sigma_1}(T_2) \includegraphics[height=1em]{crossed_arrows.pdf}
\end{align*}

where 
\begin{dmath*}
    \Delta_{id}(T_2) = - 2 v^{5} z^{5} - 17 v^{5} z^{3} - 11 v^{5} z - \frac{v^{5}}{z} + 4 v^{3} z^{9} + 40 v^{3} z^{7} + 133 v^{3} z^{5} + 150 v^{3} z^{3} + 59 v^{3} z + \frac{7 v^{3}}{z} - 5 v z^{13} - 58 v z^{11} - 255 v z^{9} - 555 v z^{7} - 689 v z^{5} - 482 v z^{3} - 161 v z - \frac{18 v}{z} + \frac{3 z^{13}}{v} + \frac{47 z^{11}}{v} + \frac{279 z^{9}}{v} + \frac{779 z^{7}}{v} + \frac{1053 z^{5}}{v} + \frac{699 z^{3}}{v} + \frac{219 z}{v} + \frac{21}{v z} + \frac{4 z^{13}}{v^{3}} + \frac{38 z^{11}}{v^{3}} + \frac{120 z^{9}}{v^{3}} + \frac{109 z^{7}}{v^{3}} - \frac{115 z^{5}}{v^{3}} - \frac{237 z^{3}}{v^{3}} - \frac{101 z}{v^{3}} - \frac{9}{v^{3} z} - \frac{3 z^{11}}{v^{5}} - \frac{35 z^{9}}{v^{5}} - \frac{144 z^{7}}{v^{5}} - \frac{254 z^{5}}{v^{5}} - \frac{177 z^{3}}{v^{5}} - \frac{38 z}{v^{5}} - \frac{2}{v^{5} z} + \frac{z^{11}}{v^{7}} + \frac{4 z^{9}}{v^{7}} + \frac{z^{7}}{v^{7}} + \frac{4 z^{5}}{v^{7}} + \frac{37 z^{3}}{v^{7}} + \frac{28 z}{v^{7}} + \frac{2}{v^{7} z} + \frac{3 z^{7}}{v^{9}} + \frac{14 z^{5}}{v^{9}} + \frac{17 z^{3}}{v^{9}} + \frac{4 z}{v^{9}}
\end{dmath*}
\begin{dmath*}
    \Delta_{\sigma_1}(T_2) = - 3 v^{5} z^{6} - 20 v^{5} z^{4} - 17 v^{5} z^{2} - 2 v^{5} + v^{3} z^{8} + 13 v^{3} z^{6} + 34 v^{3} z^{4} + 20 v^{3} z^{2} + v^{3} + 5 v z^{12} + 50 v z^{10} + 172 v z^{8} + 252 v z^{6} + 219 v z^{4} + 126 v z^{2} + 26 v - \frac{2 z^{12}}{v} - \frac{31 z^{10}}{v} - \frac{175 z^{8}}{v} - \frac{448 z^{6}}{v} - \frac{533 z^{4}}{v} - \frac{289 z^{2}}{v} - \frac{62}{v} - \frac{5 z^{12}}{v^{3}} - \frac{44 z^{10}}{v^{3}} - \frac{126 z^{8}}{v^{3}} - \frac{101 z^{6}}{v^{3}} + \frac{81 z^{4}}{v^{3}} + \frac{135 z^{2}}{v^{3}} + \frac{42}{v^{3}} + \frac{z^{12}}{v^{5}} + \frac{11 z^{10}}{v^{5}} + \frac{51 z^{8}}{v^{5}} + \frac{123 z^{6}}{v^{5}} + \frac{152 z^{4}}{v^{5}} + \frac{69 z^{2}}{v^{5}} + \frac{5}{v^{5}} + \frac{z^{12}}{v^{7}} + \frac{9 z^{10}}{v^{7}} + \frac{31 z^{8}}{v^{7}} + \frac{46 z^{6}}{v^{7}} + \frac{15 z^{4}}{v^{7}} - \frac{22 z^{2}}{v^{7}} - \frac{10}{v^{7}} - \frac{z^{8}}{v^{9}} - \frac{9 z^{6}}{v^{9}} - \frac{21 z^{4}}{v^{9}} - \frac{12 z^{2}}{v^{9}}
\end{dmath*}

Recall that we have the equations for the two variables $x_{id}, x_{\sigma_{1}} \in \mathbb{Z}[v^{\pm 1}, z^{\pm 1}]$:

\[
\begin{cases}
-vx_{id} + zv^{-1}x_{\sigma_1} = \Delta_{id}(T_2) \\
zv^{-1}x_{id} + (-v + z^{2}v^{-1})x_{\sigma_1} = \Delta_{\sigma_1}(T_2)
\end{cases}
\]

However, by solving the equations we obtain: 
\begin{dmath*}
x_{id} = \frac{1}{v^{8} z \left(v^{4} - v^{2} z^{2} - z^{2}\right)}
(2 v^{16} z^{6} + 17 v^{16} z^{4} + 11 v^{16} z^{2} + v^{16} - 4 v^{14} z^{10} - 39 v^{14} z^{8} - 130 v^{14} z^{6} - 144 v^{14} z^{4} - 58 v^{14} z^{2} - 7 v^{14} + 5 v^{12} z^{14} + 62 v^{12} z^{12} + 294 v^{12} z^{10} + 675 v^{12} z^{8} + 805 v^{12} z^{6} + 521 v^{12} z^{4} + 167 v^{12} z^{2} + 18 v^{12} - 5 v^{10} z^{16} - 66 v^{10} z^{14} - 352 v^{10} z^{12} - 1006 v^{10} z^{10} - 1720 v^{10} z^{8} - 1754 v^{10} z^{6} - 986 v^{10} z^{4} - 263 v^{10} z^{2} - 21 v^{10} + 3 v^{8} z^{16} + 45 v^{8} z^{14} + 272 v^{8} z^{12} + 834 v^{8} z^{10} + 1392 v^{8} z^{8} + 1347 v^{8} z^{6} + 745 v^{8} z^{4} + 184 v^{8} z^{2} + 9 v^{8} + 4 v^{6} z^{16} + 43 v^{6} z^{14} + 167 v^{6} z^{12} + 270 v^{6} z^{10} + 130 v^{6} z^{8} - 64 v^{6} z^{6} - 59 v^{6} z^{4} - 13 v^{6} z^{2} + 2 v^{6} - 4 v^{4} z^{14} - 47 v^{4} z^{12} - 199 v^{4} z^{10} - 378 v^{4} z^{8} - 333 v^{4} z^{6} - 144 v^{4} z^{4} - 35 v^{4} z^{2} - 2 v^{4} - 5 v^{2} z^{12} - 30 v^{2} z^{10} - 45 v^{2} z^{8} + 8 v^{2} z^{6} + 33 v^{2} z^{4} + 8 v^{2} z^{2} + 4 z^{10} + 23 z^{8} + 38 z^{6} + 16 z^{4})
\end{dmath*}

\begin{dmath*}
x_{\sigma_1} = \frac{1}{v^{8} \left(v^{4} - v^{2} z^{2} - z^{2}\right)}
(3 v^{16} z^{6} + 20 v^{16} z^{4} + 17 v^{16} z^{2} + 2 v^{16} - v^{14} z^{8} - 11 v^{14} z^{6} - 17 v^{14} z^{4} - 9 v^{14} z^{2} - 5 v^{12} z^{12} - 54 v^{12} z^{10} - 212 v^{12} z^{8} - 385 v^{12} z^{6} - 369 v^{12} z^{4} - 185 v^{12} z^{2} - 33 v^{12} + 5 v^{10} z^{14} + 60 v^{10} z^{12} + 286 v^{10} z^{10} + 730 v^{10} z^{8} + 1137 v^{10} z^{6} + 1015 v^{10} z^{4} + 450 v^{10} z^{2} + 80 v^{10} - 3 v^{8} z^{14} - 42 v^{8} z^{12} - 235 v^{8} z^{10} - 653 v^{8} z^{8} - 952 v^{8} z^{6} - 780 v^{8} z^{4} - 354 v^{8} z^{2} - 63 v^{8} - 4 v^{6} z^{14} - 39 v^{6} z^{12} - 131 v^{6} z^{10} - 160 v^{6} z^{8} - 8 v^{6} z^{6} + 85 v^{6} z^{4} + 32 v^{6} z^{2} + 4 v^{6} + 2 v^{4} z^{12} + 26 v^{4} z^{10} + 113 v^{4} z^{8} + 208 v^{4} z^{6} + 162 v^{4} z^{4} + 60 v^{4} z^{2} + 12 v^{4} - v^{2} z^{12} - 4 v^{2} z^{10} + 5 v^{2} z^{6} - 16 v^{2} z^{4} - 16 v^{2} z^{2} - 2 v^{2} - 3 z^{8} - 14 z^{6} - 17 z^{4} - 4 z^{2})
\end{dmath*}

The solutions are not in $\mathbb{Z}[v^{\pm 1}, z^{\pm 1}]$. 
This means that there is no regular isotopy between the given knot diagrams $6_1$ and $6_2$ here. 
Since they have the same writhe and Whitney index, 
they are not isotopic.

If we put $v = 1$ in $x_{id}$ and $x_{\sigma_1}$, we get 
\[
x_{id}\bigg|_{v=1} = \frac{- 2 z^{15} - 23 z^{13} - 97 z^{11} - 163 z^{9} - 38 z^{7} + 81 z^{5} + z^{3} - z}{2 z^{2} - 1}
\]
\[
x_{\sigma_1}\bigg|_{v=1} = \frac{z^{2} \left(2 z^{12} + 25 z^{10} + 112 z^{8} + 186 z^{6} + 17 z^{4} - 83 z^{2} + 9\right)}{2 z^{2} - 1}
\]
They are still not in $\mathbb{Z}[z^{\pm 1}]$. 
This means that the Alexander-Conway version of $R_{reg}^{1}$ 
can still distinguish the knot $6_1$ and $6_2$. 
\begin{remark}
    Alexander-Conway version can be calculated in polynomial time. 
\end{remark}

\subsection{More examples}
The followings are the results of the 2-cabled $R_{reg}^{1}$ on more examples. 
One could also find the computer program results at \cite{ZBT99}. 

\noindent
\begin{minipage}{\linewidth}
  \centering
  \resizebox{\linewidth}{!}{%
  \begin{tabular}{|c|c|c|}
    \hline
    & \multicolumn{2}{c|}{$R_{reg}^{1}$ distinguishes or not} \\ \hline
    & HOMFLYPT version  & Alexander-Conway version \\ \hline
    $3_1^{+}$ and $4_1$ & yes & yes \\ \hline
    $3_1^{-}$ and $3_1^{+}$ & yes & no \\ \hline
    $9_{42}$ and its mirror image & yes & no \\ \hline
    Conway knot and Kinoshita-Terasaka knot & no & no \\ \hline
  \end{tabular}%
  }
  \captionof{table}{Results of $R_{reg}^{1}$ of 2-cables}
  \label{tab:Rreg-results}
\end{minipage}
In addition, we provide an example illustrating the information 
that quantum equations yield about the regular isotopies 
connecting the standard knot diagram of $4_1$ and its mirror image.

\begin{itemize}
    \begin{item}
        $3_1^{+}$ and $4_1$ 

        \noindent
        \begin{minipage}{\linewidth}
        \captionsetup{type=figure} 
        \centering
        \begin{subfigure}{0.4\linewidth}
            \includegraphics[width=\linewidth]{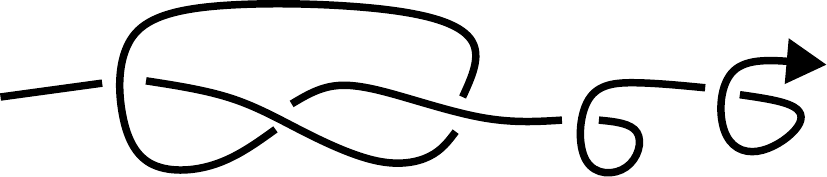}
            \caption{$3_1^{+}$}
        \end{subfigure}\hfill
        \begin{subfigure}{0.5\linewidth}
            \includegraphics[width=\linewidth]{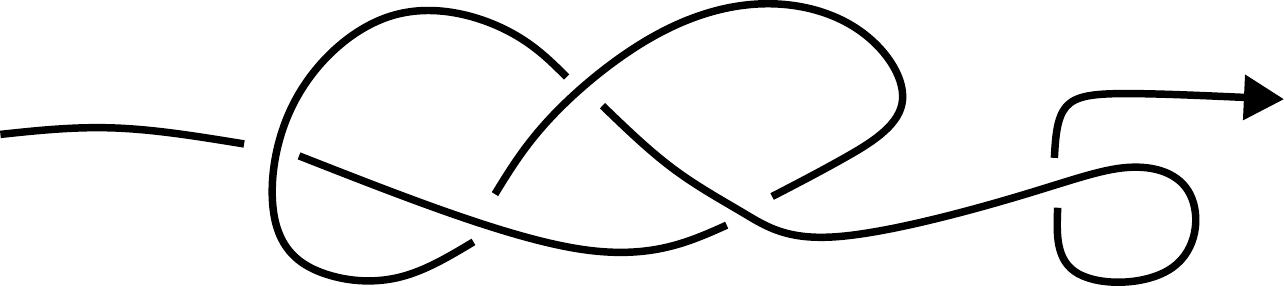}
            \caption{$4_1$}
        \end{subfigure}
        \captionof{figure}{$3_1^{+}$ and $4_1$ with writhe $=1$ and Whitney index $=1$}
        \label{fig:trefoil_figure8}
        \end{minipage}
        \par
        From pushing the tangle $T_2$ we get: 
        \begin{dmath*}
            x_{id} = \frac{1}{v^{6} z \left(- v^{4} + v^{2} z^{2} + z^{2}\right)}\cdot(- 5 v^{16} z^{4} - 12 v^{16} z^{2} - 2 v^{16} + 12 v^{14} z^{6} + 46 v^{14} z^{4} + 37 v^{14} z^{2} + 6 v^{14} - 8 v^{12} z^{8} - 41 v^{12} z^{6} - 60 v^{12} z^{4} - 33 v^{12} z^{2} - 9 v^{12} + v^{10} z^{10} + 4 v^{10} z^{8} - 2 v^{10} z^{6} - 8 v^{10} z^{4} + 21 v^{10} z^{2} + 8 v^{10} + v^{8} z^{10} + 9 v^{8} z^{8} + 16 v^{8} z^{6} - 8 v^{8} z^{4} - v^{8} z^{2} - 2 v^{8} + 3 v^{6} z^{10} + 18 v^{6} z^{8} + 19 v^{6} z^{6} - 16 v^{6} z^{4} - 12 v^{6} z^{2} - 2 v^{6} + 3 v^{4} z^{10} + 15 v^{4} z^{8} + 14 v^{4} z^{6} - 2 v^{4} z^{4} + 2 v^{4} z^{2} + v^{4} + 4 v^{2} z^{6} + 9 v^{2} z^{4} + 7 v^{2} z^{2} - 3 z^{6} - 8 z^{4} - 2 z^{2})
        \end{dmath*}
        \begin{dmath*}
            x_{\sigma_1} = \frac{1}{v^{6} \left(- v^{4} + v^{2} z^{2} + z^{2}\right)}(- 5 v^{16} z^{4} - 17 v^{16} z^{2} - 10 v^{16} + 12 v^{14} z^{6} + 59 v^{14} z^{4} + 77 v^{14} z^{2} + 23 v^{14} - 8 v^{12} z^{8} - 51 v^{12} z^{6} - 100 v^{12} z^{4} - 62 v^{12} z^{2} - 19 v^{12} + v^{10} z^{10} + 6 v^{10} z^{8} + 3 v^{10} z^{6} - 21 v^{10} z^{4} + 6 v^{10} z^{2} + 4 v^{10} + v^{8} z^{10} + 13 v^{8} z^{8} + 44 v^{8} z^{6} + 42 v^{8} z^{4} + 23 v^{8} z^{2} + 2 v^{8} - v^{6} z^{8} - 10 v^{6} z^{6} - 23 v^{6} z^{4} - 10 v^{6} z^{2} - v^{4} z^{8} - 5 v^{4} z^{6} - 10 v^{4} z^{4} - 4 v^{4} z^{2} + v^{4} - 6 v^{2} z^{4} - 5 v^{2} z^{2} - v^{2} + 3 z^{4} + 5 z^{2})
        \end{dmath*}
        if we put $v = 1$, we get
        \begin{dmath*}
            x_{id}\bigg|_{v=1} = \frac{z \left(8 z^{8} + 38 z^{6} + 19 z^{4} - 52 z^{2} + 7\right)}{2 z^{2} - 1}
        \end{dmath*}
        \begin{dmath*}
            x_{\sigma_1}\bigg|_{v=1} = \frac{z^{2} \left(2 z^{8} + 9 z^{6} - 7 z^{4} - 61 z^{2} + 13\right)}{2 z^{2} - 1}
        \end{dmath*}
        The results mean that both Alexander-Conway version and HOMFLYPT version can distinguish the knots $3_1^{+}$ and $4_1$. 
    \end{item}

    \begin{item}
        $3_1^{-}$ and $3_1^{+}$ 

        \noindent
        \begin{minipage}{\linewidth}
        \captionsetup{type=figure} 
        \centering
        \begin{subfigure}{0.4\linewidth}
            \includegraphics[width=\linewidth]{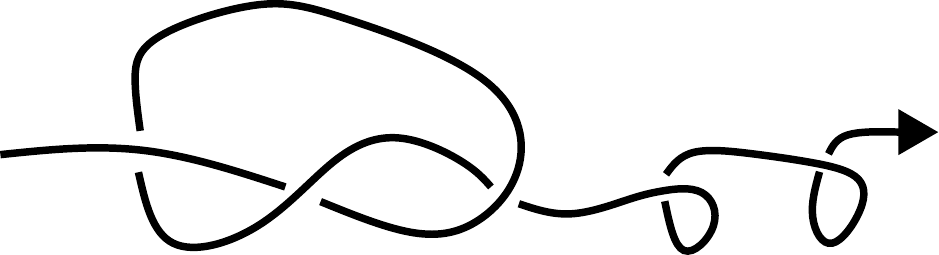}
            \caption{$3_1^{-}$}
        \end{subfigure}\hfill
        \begin{subfigure}{0.5\linewidth}
            \includegraphics[width=\linewidth]{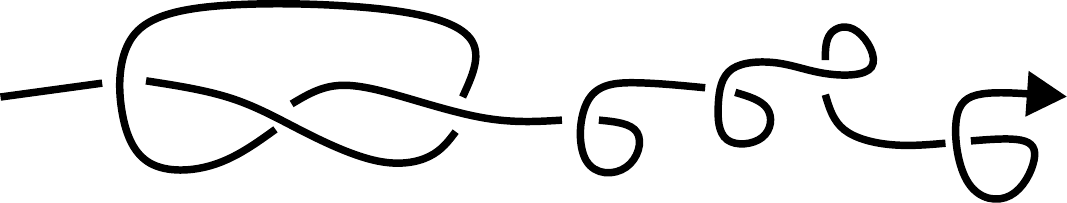}
            \caption{$3_1^{+}$}
        \end{subfigure}
        \captionof{figure}{$3_1^{-}$ and $3_1^{+}$ with writhe $=-1$ and Whitney index $=1$}
        \label{fig:trefoil_left_right}
        \end{minipage}
        \par
        From pushing the tangle $T_2$ we get 
        \begin{dmath*}
            x_{id} = \frac{1}{v^{10} z \left(- v^{4} + v^{2} z^{2} + z^{2}\right)}(6 v^{20} z^{2} + 4 v^{20} - 2 v^{18} z^{6} - 14 v^{18} z^{4} - 27 v^{18} z^{2} - 16 v^{18} + 2 v^{16} z^{8} + 12 v^{16} z^{6} + 41 v^{16} z^{4} + 51 v^{16} z^{2} + 20 v^{16} - 2 v^{14} z^{8} - 20 v^{14} z^{6} - 34 v^{14} z^{4} - 25 v^{14} z^{2} - 8 v^{14} - v^{12} z^{8} - 8 v^{12} z^{6} - 15 v^{12} z^{4} - 8 v^{12} z^{2} + 3 v^{10} z^{8} + 19 v^{10} z^{6} + 34 v^{10} z^{4} + 23 v^{10} z^{2} + 4 v^{10} + 3 v^{8} z^{12} + 23 v^{8} z^{10} + 51 v^{8} z^{8} + 9 v^{8} z^{6} - 81 v^{8} z^{4} - 64 v^{8} z^{2} - 10 v^{8} - 3 v^{6} z^{14} - 24 v^{6} z^{12} - 56 v^{6} z^{10} - v^{6} z^{8} + 160 v^{6} z^{6} + 202 v^{6} z^{4} + 88 v^{6} z^{2} + 8 v^{6} - 3 v^{4} z^{14} - 27 v^{4} z^{12} - 91 v^{4} z^{10} - 156 v^{4} z^{8} - 165 v^{4} z^{6} - 124 v^{4} z^{4} - 43 v^{4} z^{2} - 2 v^{4} - v^{2} z^{12} - 8 v^{2} z^{10} - 21 v^{2} z^{8} - 19 v^{2} z^{6} - 5 v^{2} z^{4} + 2 v^{2} z^{2} + 4 z^{10} + 24 z^{8} + 41 z^{6} + 19 z^{4} + z^{2})
        \end{dmath*}
        \begin{dmath*}
            x_{\sigma_1} = \frac{1}{v^{10} \left(- v^{4} + v^{2} z^{2} + z^{2}\right)}(6 v^{20} z^{2} + 10 v^{20} - 2 v^{18} z^{6} - 18 v^{18} z^{4} - 42 v^{18} z^{2} - 23 v^{18} + 4 v^{16} z^{8} + 25 v^{16} z^{6} + 50 v^{16} z^{4} + 30 v^{16} z^{2} + 11 v^{16} - 2 v^{14} z^{10} - 11 v^{14} z^{8} - 13 v^{14} z^{6} + 13 v^{14} z^{4} + 15 v^{14} z^{2} + 2 v^{14} - 2 v^{12} z^{10} - 13 v^{12} z^{8} - 26 v^{12} z^{6} - 23 v^{12} z^{4} - 12 v^{12} z^{2} - v^{10} z^{6} - 6 v^{10} z^{4} - 8 v^{10} z^{2} - 5 v^{10} - 3 v^{8} z^{10} - 20 v^{8} z^{8} - 34 v^{8} z^{6} + 9 v^{8} z^{4} + 58 v^{8} z^{2} + 25 v^{8} + 3 v^{6} z^{12} + 21 v^{6} z^{10} + 38 v^{6} z^{8} - 21 v^{6} z^{6} - 128 v^{6} z^{4} - 115 v^{6} z^{2} - 33 v^{6} + 3 v^{4} z^{12} + 24 v^{4} z^{10} + 70 v^{4} z^{8} + 105 v^{4} z^{6} + 98 v^{4} z^{4} + 65 v^{4} z^{2} + 14 v^{4} + v^{2} z^{10} + 7 v^{2} z^{8} + 15 v^{2} z^{6} + 9 v^{2} z^{4} + v^{2} z^{2} - v^{2} - 4 z^{8} - 20 z^{6} - 25 z^{4} - 6 z^{2})
        \end{dmath*}
        If we put $v = 1$, we get
        \begin{dmath*}
            x_{id}\bigg|_{v=1} = z \left(- 3 z^{10} - 26 z^{8} - 77 z^{6} - 89 z^{4} - 31 z^{2} - 4\right)
        \end{dmath*}
        \begin{dmath*}
            x_{\sigma_1}\bigg|_{v=1} = z^{2} \left(3 z^{8} + 21 z^{6} + 46 z^{4} + 37 z^{2} + 8\right)
        \end{dmath*}
        $x_{id}\bigg|_{v=1}$ and $x_{\sigma_1}\bigg|_{v=1}$ are the same as the results from pushing the tangle $T_1$ (further calculations required).
        The results mean that the Alexander-Conway version of $R_{reg}^{1}$ does not detect the chirality of trefoil but the HOMFLYPT version does. 
    \end{item}

    \begin{item}
        $9_{42}$ and its mirror image ($9_{42}$ is chiral)

        \noindent
        \begin{minipage}{\linewidth}
        \captionsetup{type=figure} 
        \centering
        \begin{subfigure}{0.4\linewidth}
            \includegraphics[width=\linewidth]{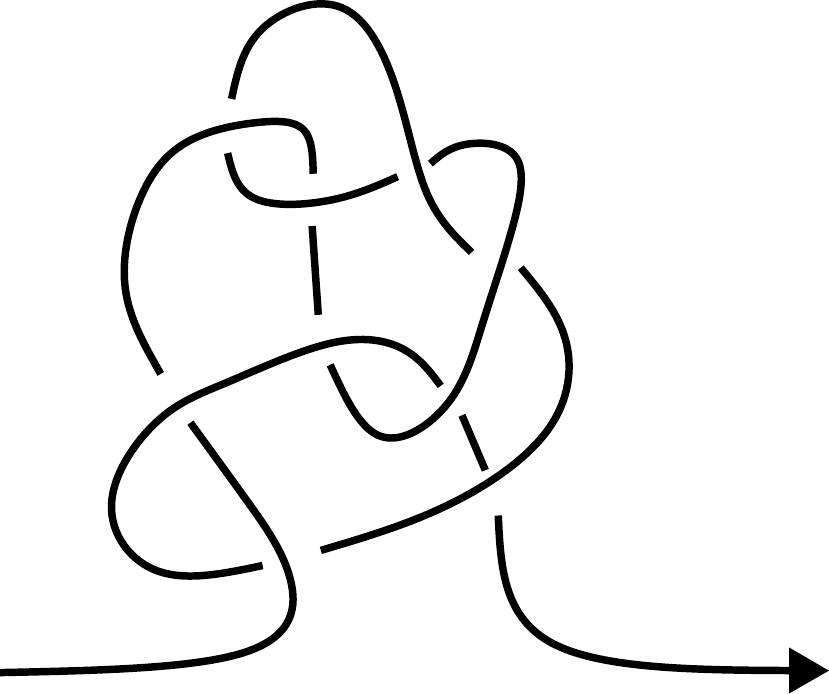}
            \caption{$9_{42}$}
        \end{subfigure}\hfill
        \begin{subfigure}{0.5\linewidth}
            \includegraphics[width=\linewidth]{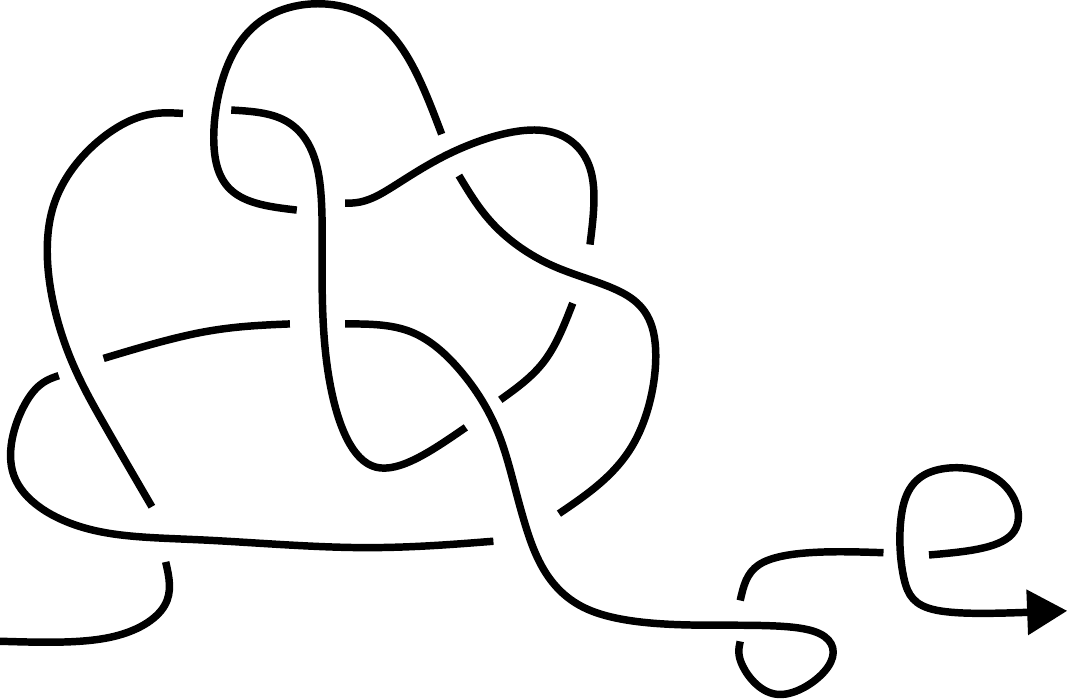}
            \caption{the mirror of $9_{42}$}
        \end{subfigure}
        \captionof{figure}{$9_{42}$ and the mirror of $9_{42}$ with writhe $=1$ and Whitney index $=1$}
        \label{fig:9_42_and_mirror}
        \end{minipage}
        \par
        From pushing the tangle $T_2$ we get 
        \begin{dmath*}
            x_{id} = \frac{1}{v^{4} z \left(- v^{4} + v^{2} z^{2} + z^{2}\right)}(- 8 v^{16} z^{8} - 50 v^{16} z^{6} - 88 v^{16} z^{4} - 39 v^{16} z^{2} + 8 v^{14} z^{12} + 99 v^{14} z^{10} + 429 v^{14} z^{8} + 795 v^{14} z^{6} + 613 v^{14} z^{4} + 149 v^{14} z^{2} - 2 v^{14} - 16 v^{12} z^{14} - 196 v^{12} z^{12} - 914 v^{12} z^{10} - 2044 v^{12} z^{8} - 2271 v^{12} z^{6} - 1150 v^{12} z^{4} - 207 v^{12} z^{2} + 9 v^{12} + 7 v^{10} z^{16} + 95 v^{10} z^{14} + 514 v^{10} z^{12} + 1407 v^{10} z^{10} + 2054 v^{10} z^{8} + 1600 v^{10} z^{6} + 687 v^{10} z^{4} + 97 v^{10} z^{2} - 17 v^{10} + v^{8} z^{18} + 12 v^{8} z^{16} + 64 v^{8} z^{14} + 223 v^{8} z^{12} + 543 v^{8} z^{10} + 770 v^{8} z^{8} + 512 v^{8} z^{6} + 177 v^{8} z^{4} + 56 v^{8} z^{2} + 17 v^{8} - 7 v^{6} z^{16} - 87 v^{6} z^{14} - 420 v^{6} z^{12} - 1002 v^{6} z^{10} - 1258 v^{6} z^{8} - 852 v^{6} z^{6} - 343 v^{6} z^{4} - 74 v^{6} z^{2} - 9 v^{6} - v^{4} z^{18} - 10 v^{4} z^{16} - 36 v^{4} z^{14} - 64 v^{4} z^{12} - 83 v^{4} z^{10} - 67 v^{4} z^{8} + 27 v^{4} z^{6} + 20 v^{4} z^{4} + 14 v^{4} z^{2} + 2 v^{4} + v^{2} z^{14} + 15 v^{2} z^{12} + 87 v^{2} z^{10} + 230 v^{2} z^{8} + 269 v^{2} z^{6} + 92 v^{2} z^{4} + 9 v^{2} z^{2} - z^{12} - 10 z^{10} - 36 z^{8} - 58 z^{6} - 39 z^{4} - 2 z^{2})
        \end{dmath*}
        \begin{dmath*}
            x_{\sigma_1} = \frac{1}{v^{4} \left(- v^{4} + v^{2} z^{2} + z^{2}\right)}(- 8 v^{16} z^{8} - 58 v^{16} z^{6} - 132 v^{16} z^{4} - 98 v^{16} z^{2} - 10 v^{16} + 8 v^{14} z^{12} + 107 v^{14} z^{10} + 524 v^{14} z^{8} + 1171 v^{14} z^{6} + 1200 v^{14} z^{4} + 475 v^{14} z^{2} + 51 v^{14} - 16 v^{12} z^{14} - 214 v^{12} z^{12} - 1119 v^{12} z^{10} - 2906 v^{12} z^{8} - 3947 v^{12} z^{6} - 2738 v^{12} z^{4} - 943 v^{12} z^{2} - 125 v^{12} + 7 v^{10} z^{16} + 102 v^{10} z^{14} + 604 v^{10} z^{12} + 1877 v^{10} z^{10} + 3349 v^{10} z^{8} + 3650 v^{10} z^{6} + 2583 v^{10} z^{4} + 1066 v^{10} z^{2} + 178 v^{10} + v^{8} z^{18} + 15 v^{8} z^{16} + 98 v^{8} z^{14} + 363 v^{8} z^{12} + 777 v^{8} z^{10} + 761 v^{8} z^{8} - 195 v^{8} z^{6} - 1023 v^{8} z^{4} - 695 v^{8} z^{2} - 151 v^{8} - 7 v^{6} z^{16} - 97 v^{6} z^{14} - 529 v^{6} z^{12} - 1433 v^{6} z^{10} - 1996 v^{6} z^{8} - 1321 v^{6} z^{6} - 255 v^{6} z^{4} + 180 v^{6} z^{2} + 71 v^{6} - v^{4} z^{18} - 13 v^{4} z^{16} - 64 v^{4} z^{14} - 133 v^{4} z^{12} - 11 v^{4} z^{10} + 469 v^{4} z^{8} + 819 v^{4} z^{6} + 483 v^{4} z^{4} + 78 v^{4} z^{2} - 14 v^{4} + 2 v^{2} z^{14} + 18 v^{2} z^{12} + 45 v^{2} z^{10} - 19 v^{2} z^{8} - 199 v^{2} z^{6} - 229 v^{2} z^{4} - 72 v^{2} z^{2} + 4 z^{8} + 24 z^{6} + 42 z^{4} + 16 z^{2})
        \end{dmath*}
        If we put $v = 1$, we get
        \begin{dmath*}
            x_{id}\bigg|_{v=1} = z \left(z^{12} + 11 z^{10} + 45 z^{8} + 86 z^{6} + 78 z^{4} + 25 z^{2} - 3\right)
        \end{dmath*}
        \begin{dmath*}
            x_{\sigma_1}\bigg|_{v=1} = z^{2} \left(z^{12} + 13 z^{10} + 65 z^{8} + 154 z^{6} + 166 z^{4} + 55 z^{2} - 7\right)
        \end{dmath*}
        $x_{id}\bigg|_{v=1}$ and $x_{\sigma_1}\bigg|_{v=1}$ are the same as the results from pushing the full twist $T_1$ (further calculations required).
        The results mean that the Alexander-Conway version of $R_{reg}^{1}$ does not detect the chirality of $9_{42}$ but the HOMFLYPT version does. 
        It should be noticed that $9_{42}$ and its mirror image share the same HOMFLYPT polynomial. 
        In \cite{MS87} Morton and Short showed that the HOMFLYPT polynomial of the 2-cable of $9_{42}$ detects the chirality of $9_{42}$. 
        
        Let us specialise $v = i$. For the full twist $T_1$, we obtain: 
        \[
        x_{id}\bigg|_{v=i} = - 2 z^{15} - 11 z^{13} + 51 z^{11} + 527 z^{9} + 1512 z^{7} + 1924 z^{5} + 1041 z^{3} + 103 z - \frac{56}{z}
        \]
        \[
        x_{\sigma_1}\bigg|_{v=i} = - 2 z^{16} - 15 z^{14} + 33 z^{12} + 685 z^{10} + 2874 z^{8} + 5794 z^{6} + 6123 z^{4} + 3163 z^{2} + 600
        \]
        For the tangle $T_2$ we obtain : 
        \[
        x_{id}\bigg|_{v=i} = - 2 z^{15} - 3 z^{13} + 155 z^{11} + 1055 z^{9} + 2840 z^{7} + 3652 z^{5} + 2129 z^{3} + 359 z - \frac{56}{z}
        \]
        \[
        x_{\sigma_1}\bigg|_{v=i} = - 2 z^{16} - 11 z^{14} + 85 z^{12} + 949 z^{10} + 3538 z^{8} + 6658 z^{6} + 6667 z^{4} + 3291 z^{2} + 600
        \]

        These variables are all Laurent polynomials but they are different. 
        Consequently, this is another proof that these knots are not isotopic.

    \end{item}

    \begin{item}
        The mutants: the Conway knot and the Kinoshita-Terasaka knot 

        \noindent
        \begin{minipage}{\linewidth}
        \captionsetup{type=figure} 
        \centering
        \begin{subfigure}{0.47\linewidth}
            \includegraphics[width=\linewidth]{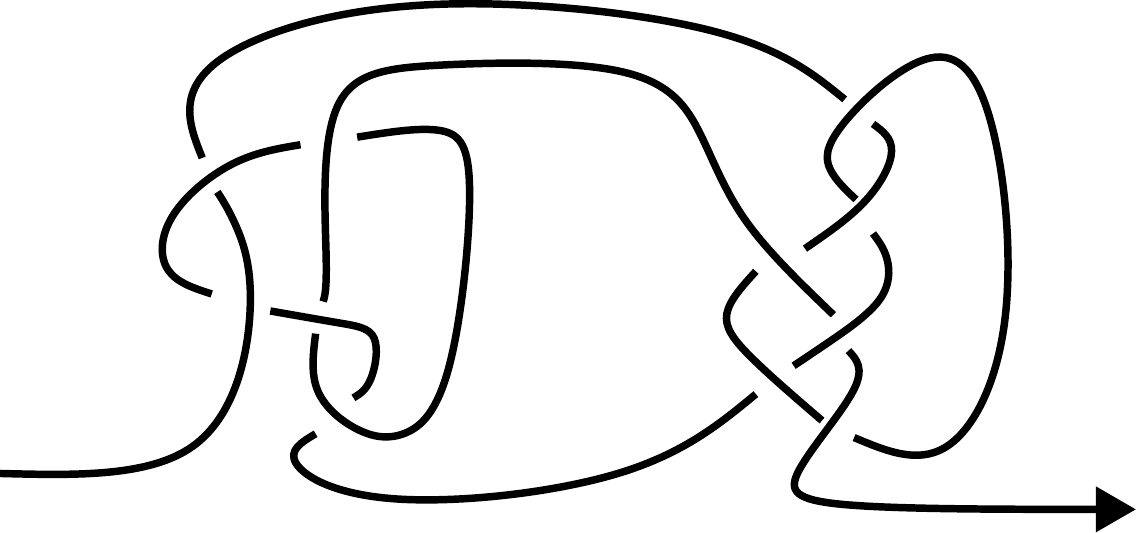}
            \caption{the Conway knot}
        \end{subfigure}\hfill
        \begin{subfigure}{0.45\linewidth}
            \includegraphics[width=\linewidth]{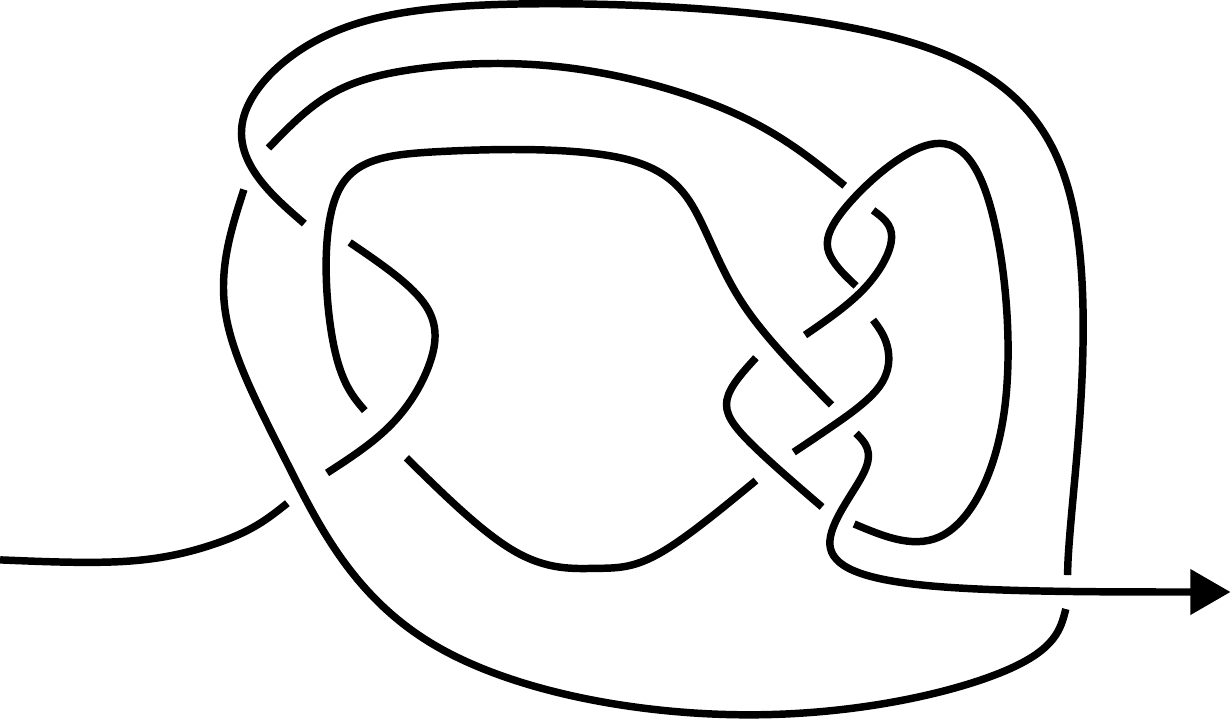}
            \caption{the Kinoshita-Terasaka knot }
        \end{subfigure}
        \captionof{figure}{The mutants with writhe $=1$ and Whitney index $=1$}
        \label{fig:mutants}
        \end{minipage}
        \par
        From pushing the tangle $T_2$ we get 
        \begin{dmath*}
            x_{id} = - 2 v^{6} z^{15} - 26 v^{6} z^{13} - 130 v^{6} z^{11} - 314 v^{6} z^{9} - 379 v^{6} z^{7} - 213 v^{6} z^{5} - 39 v^{6} z^{3} + 7 v^{6} z + \frac{2 v^{6}}{z} + v^{4} z^{19} + 19 v^{4} z^{17} + 152 v^{4} z^{15} + 665 v^{4} z^{13} + 1731 v^{4} z^{11} + 2735 v^{4} z^{9} + 2565 v^{4} z^{7} + 1322 v^{4} z^{5} + 288 v^{4} z^{3} - 14 v^{4} z - \frac{11 v^{4}}{z} - v^{2} z^{21} - 21 v^{2} z^{19} - 191 v^{2} z^{17} - 984 v^{2} z^{15} - 3153 v^{2} z^{13} - 6500 v^{2} z^{11} - 8630 v^{2} z^{9} - 7183 v^{2} z^{7} - 3491 v^{2} z^{5} - 805 v^{2} z^{3} - 6 v^{2} z + \frac{22 v^{2}}{z} + z^{23} + 22 z^{21} + 210 z^{19} + 1146 z^{17} + 3972 z^{15} + 9195 z^{13} + 14536 z^{11} + 15683 z^{9} + 11232 z^{7} + 4958 z^{5} + 1116 z^{3} + 39 z - \frac{21}{z} - \frac{z^{23}}{v^{2}} - \frac{21 z^{21}}{v^{2}} - \frac{191 z^{19}}{v^{2}} - \frac{991 z^{17}}{v^{2}} - \frac{3258 z^{15}}{v^{2}} - \frac{7141 z^{13}}{v^{2}} - \frac{10693 z^{11}}{v^{2}} - \frac{10983 z^{9}}{v^{2}} - \frac{7569 z^{7}}{v^{2}} - \frac{3265 z^{5}}{v^{2}} - \frac{730 z^{3}}{v^{2}} - \frac{33 z}{v^{2}} + \frac{10}{v^{2} z} + \frac{z^{19}}{v^{4}} + \frac{17 z^{17}}{v^{4}} + \frac{120 z^{15}}{v^{4}} + \frac{460 z^{13}}{v^{4}} + \frac{1057 z^{11}}{v^{4}} + \frac{1516 z^{9}}{v^{4}} + \frac{1349 z^{7}}{v^{4}} + \frac{701 z^{5}}{v^{4}} + \frac{182 z^{3}}{v^{4}} + \frac{13 z}{v^{4}} - \frac{2}{v^{4} z} - \frac{z^{7}}{v^{6}} - \frac{4 z^{5}}{v^{6}} - \frac{z^{3}}{v^{6}} - \frac{z}{v^{6}}
        \end{dmath*}
        \begin{dmath*}
            x_{\sigma_1} = - v^{6} z^{16} - 12 v^{6} z^{14} - 50 v^{6} z^{12} - 69 v^{6} z^{10} + 67 v^{6} z^{8} + 287 v^{6} z^{6} + 297 v^{6} z^{4} + 129 v^{6} z^{2} + 20 v^{6} - 2 v^{4} z^{16} - 35 v^{4} z^{14} - 250 v^{4} z^{12} - 942 v^{4} z^{10} - 2026 v^{4} z^{8} - 2532 v^{4} z^{6} - 1809 v^{4} z^{4} - 685 v^{4} z^{2} - 106 v^{4} - v^{2} z^{22} - 20 v^{2} z^{20} - 169 v^{2} z^{18} - 780 v^{2} z^{16} - 2106 v^{2} z^{14} - 3210 v^{2} z^{12} - 1999 v^{2} z^{10} + 1569 v^{2} z^{8} + 4078 v^{2} z^{6} + 3370 v^{2} z^{4} + 1336 v^{2} z^{2} + 213 v^{2} + z^{24} + 23 z^{22} + 230 z^{20} + 1312 z^{18} + 4705 z^{16} + 11000 z^{14} + 16727 z^{12} + 15710 z^{10} + 7464 z^{8} - 445 z^{6} - 2589 z^{4} - 1286 z^{2} - 221 - \frac{z^{24}}{v^{2}} - \frac{22 z^{22}}{v^{2}} - \frac{211 z^{20}}{v^{2}} - \frac{1161 z^{18}}{v^{2}} - \frac{4057 z^{16}}{v^{2}} - \frac{9395 z^{14}}{v^{2}} - \frac{14524 z^{12}}{v^{2}} - \frac{14542 z^{10}}{v^{2}} - \frac{8472 z^{8}}{v^{2}} - \frac{1702 z^{6}}{v^{2}} + \frac{1023 z^{4}}{v^{2}} + \frac{709 z^{2}}{v^{2}} + \frac{129}{v^{2}} + \frac{z^{20}}{v^{4}} + \frac{18 z^{18}}{v^{4}} + \frac{135 z^{16}}{v^{4}} + \frac{548 z^{14}}{v^{4}} + \frac{1307 z^{12}}{v^{4}} + \frac{1840 z^{10}}{v^{4}} + \frac{1384 z^{8}}{v^{4}} + \frac{282 z^{6}}{v^{4}} - \frac{323 z^{4}}{v^{4}} - \frac{226 z^{2}}{v^{4}} - \frac{41}{v^{4}} + \frac{z^{12}}{v^{6}} + \frac{10 z^{10}}{v^{6}} + \frac{36 z^{8}}{v^{6}} + \frac{63 z^{6}}{v^{6}} + \frac{67 z^{4}}{v^{6}} + \frac{35 z^{2}}{v^{6}} + \frac{6}{v^{6}}
        \end{dmath*}
    
        $x_{id}$ and $x_{\sigma_1}$ are the same as the results from pushing the tangle $T_1$ (further calculations required).
        This means that the HOMFLYPT version of $R_{reg}^{1}$ does not distinguish the mutants: the Conway knot and the Kinoshita-Terasaka knot.
        In \cite{MS87} Morton and Short showed that the HOMFLYPT polynomial of the 2-cable of $9_{42}$ cannot distinguish them either.  
    \end{item}

    \begin{item}
        $4_1$ and its mirror image

        \noindent
        \begin{minipage}{\linewidth}
        \captionsetup{type=figure} 
        \centering
        \begin{subfigure}{0.4\linewidth}
            \includegraphics[width=\linewidth]{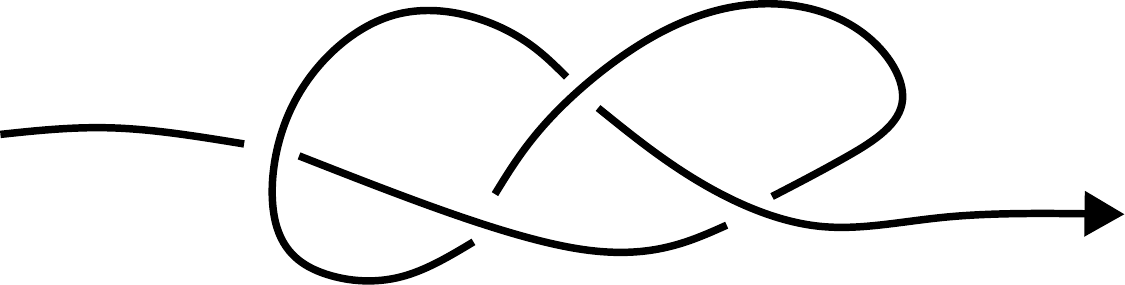}
            \caption{$4_1$}
        \end{subfigure}\hfill
        \begin{subfigure}{0.4\linewidth}
            \includegraphics[width=\linewidth]{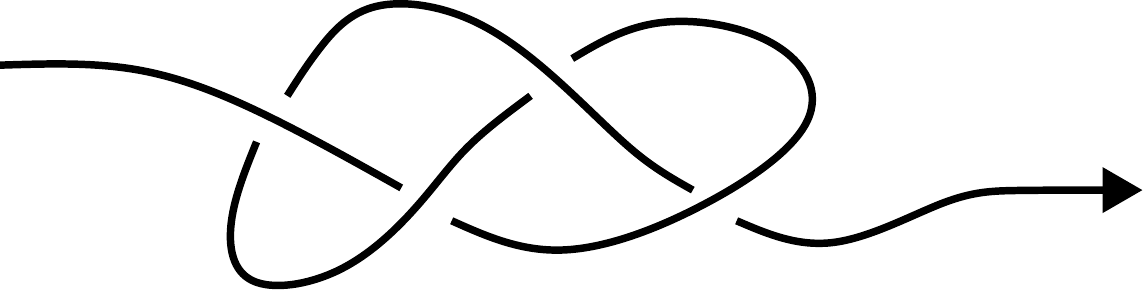}
            \caption{the mirror of $4_1$}
        \end{subfigure}
        \captionof{figure}{$4_{1}$ and the mirror of $4_{1}$ with writhe $=$ Whitney index $=0$}
        \label{fig:4_1_and_mirror}
        \end{minipage}
        \par
        It is well known that the figure-eight knot $4_1$ is isotopic to its mirror image. 
        For each regular isotopy which connects the standard diagram $4_1$ with its mirror image, the quantum equations for using $T_1$ or $T_2$ give now: 
        \[
        x_{id} = v^{2} z - z^{5} - 4 z^{3} - \frac{z^{5}}{v^{2}} - \frac{4 z^{3}}{v^{2}} + \frac{z}{v^{4}}
        \]
        \[
        x_{\sigma_1} = - v^{4} + v^{2} z^{4} + 5 v^{2} z^{2} + v^{2} - z^{6} - 4 z^{4} - \frac{z^{6}}{v^{2}} - \frac{5 z^{4}}{v^{2}} - \frac{4 z^{2}}{v^{2}} - \frac{1}{v^{2}} + \frac{z^{2}}{v^{4}} + \frac{1}{v^{4}}
        \]
        Notice that we know this without knowing the regular isotopy. 

    \end{item}
\end{itemize}

Institut de Math\'ematiques de Toulouse, UMR 5219

Universit\'e de Toulouse

118, route de Narbonne 

31062 Toulouse Cedex 09, France

thomas.fiedler@math.univ-toulouse.fr

butian.zhang@math.univ-toulouse.fr
\end{document}